\documentclass[12pt,twoside]{article}
\usepackage{amsmath,amssymb}
\newcommand{\N}{\mathbb{N}}
\newcommand{\Z}{\mathbb{Z}}
\newcommand{\Q}{\mathbb{Q}}

\newcommand{\C}{\mathbb{C}}
\newcommand{\F}{\mathbb{F}}
\renewcommand{\L}{\mathbb{L}}
\renewcommand{\H}{\mathbb{H}}
\renewcommand{\P}{\mathbb{P}}
\newcommand{\lra}{\longrightarrow}
\usepackage[arrow,matrix,curve]{xy}
\newtheorem{theorem}{Theorem}[section]

\newtheorem{proposition}{Proposition}[section]

\newtheorem{definition}{Definition}[section]

\newtheorem{conjecture}{Conjecture}[section]
\begin{document}

\title{Frobenius polynomials for Calabi-Yau equations}

\author{Kira Samol and Duco van Straten}

\date{15 September 2008}

\pagestyle{myheadings}
\markboth{K.~Samol and D.~van Straten}%
{Frobenius Polynomials}

\maketitle

\begin{abstract}
We describe a variation of Dwork' s unit-root method to
determine the degree four Frobenius polynomial
for members of a 1-modulus Calabi-Yau family over $\P^1$ 
in terms of the holomorphic period near a point of maximal
unipotent monodromy.
The method is illustrated on a couple of examples from the
list \cite{AESZ}. 
For singular points we find that the Frobenius polynomial 
splits in a product of two linear factors and 
a quadratic part $1-a_pT+p^3T^2$. We identify weight four 
modular forms which reproduce the $a_p$ as Fourier coefficients.
\end{abstract}

\bigskip

\section{Introduction}

Given a projective morphism $f:X \lra \P^1$ with smooth generic $n-1$-dimensional fibre, the sheaf
$R^{n-1}f_*(\Q_X)$ restricts to a $\Q$-local system $\H$ over the smooth locus $S \subset \P^1$ of $f$ and
hence determines, after the choice of a base-point $s_0 \in S$, a monodromy representation 
$\pi_1(S,s_0) \lra Aut(\H_{s_0})$. The local system $\H$ carries a non-degenerate $(-1)^{n-1}$-symmetric pairing 
\[<-,->: \H \otimes \H \lra \Q_S\] induced by the intersection form in the fibres. Hence we can identify $\H$ with
its dual and the monodromy representation lands in a symplectic ($n-1$ odd) or orthogonal group ($n-1$ even). The primitive part of
$\H$ underlies a variation of Hodge structures (VHS), polarised by $<-,->$, see \cite{peterssteenbrink}. 

We call a sub-VHS $\mathbb{L} \subset \H$ a {\em $CY(n)$-local system} if the local monodromy around 
$0\in \P^1 \setminus S$ is unipotent and consists of a single Jordan block of size $n$. Hence, $\L$ is 
irreducible of rank $n$ and the non-vanishing sub-quotients $Gr^W_{2k}$ ($k=0,1,\ldots,n-1$) of the monodromy 
weight filtration all have dimension equal to one. The Hodge filtration $F^{\cdot}$ of the limiting mixed Hodge
structure at $0$ is opposite to the weight filtration \cite{deligne},\cite{morrison}. If $\omega$ is a section of the smallest Hodge space
$F^{n-1}$ and $\gamma$ a local section of $\mathbb{L}$ near $0$, then the {\em period function} 
\[f_0:=<\gamma,\omega>\]
is holomorphic near $0$ and satisfies a linear differential equation of order $n$, called the associated {\em Picard-Fuchs} equation.\\  
We call a linear differential operator of order $n$
\[P =\frac{d}{dx}^{n}+a_{n-1}(x)\frac{d}{dx}^{n-1}+\ldots+a_0(x) \in \mathbb{Q}(x)[\frac{d}{dx}]\] a 
{\em CY(n)-operator}  if
\begin{enumerate}
\item $P$ has maximal unipotent monodromy at $0$ (MUM).
\item $P$ is self-dual in the sense that
\[P =(-1)^n \exp(-\frac{2}{n}\int a_{n-1})\circ P^* \circ \exp(\frac{2}{n}\int a_{n-1}),\]
where $\circ$ means the composition of differential operators.
\item $P$ has a convergent power series solution $f_0(x) \in \Z[[x]]$ with $f_0(0)=1$.
\end{enumerate}

The first condition implies that the operator $P$ is irreducible and can (after left multiplication by $x^n$) 
be written in the form 
\[\theta^n+xP_1(\theta)+x^2P_2(\theta)+\ldots+x^dP_d(\theta),\] 
where $\theta:=x\frac{d}{dx}$ and $P_i(\theta) \in \mathbb{Q}[\theta]$ is a polynomial in $\theta$. 
We remark that $\exp(\int a_{n-1}) \in \mathbb{Q}(x)$ precisely if the differential Galois group of ${\cal P}$ belongs to $SL(n)$. In the second condition $P^*$ is the formal adjoint of $P$. The condition is equivalent to the condition that
the transformed operator
\[\tilde{P}=\exp(\frac{1}{n}\int a_{n-1})\circ P \circ \exp(-\frac{1}{n}\int a_{n-1})=\frac{d}{dx}^n+0\frac{d}{dx}^{n-1}+\ldots\]
satisfies
\[\tilde{P}=(-1)^n \tilde{P}^*\]
which translates into $\lfloor (n-1)/2 \rfloor$ differential-polynomial conditions on the coefficients $a_i$.
These express the conditions that the differential Galois group of $P$ is in the symplectic or orthogonal group. For $n=4$ one finds the condition of \cite{almkvistzudilin}:
\[a_1=\frac{1}{2}a_2a_3-\frac{1}{8}a_3^3+a_2'-\frac{3}{4}a_3a_3-\frac{1}{2}a_3''\]
If the Yukawa coupling is non-constant, then the differential Galois group is $Gal(P)^{0}=Sp(4)$ in general, see \cite{bogner}.
In \cite{AESZ} one finds a list with more than $350$ examples of such fourth order operators.

Because of the MUM-condition, the solution $f_0(x)$ from the third condition is unique and conversely determines the operator $P$. As $f_0$ is a $G$-function, the operator $P$ is a $G$-operator and hence by a theorem of Katz is of fuchsian type with rational exponents, see \cite{andre}.

A Picard-Fuchs operator that arises from a geometrical situation as sketched above will satisfy the first two conditions and the period function $f_0$ will be a so called  G-function, see \cite{andre}.
It would therefore perhaps seem  more  natural to require $f_0$ to be a $G$-function. However, requiring integrality of the solution covers all interesting examples and helps fixing the coordinate $x$.  
In \cite{almkvistzudilin} for $n=4$ further integrality properties for the  mirror map and Yukawa coupling were required. 

CY(2)-operators arise from families of elliptic curves, CY(3)-operators arise from families of $K3$ with Picard-number $19$ with a point of maximal degeneration (type III in the terminology of \cite{friedmanmorrison}).
CY(4)-operators arise from families of Calabi-Yau threefolds with $h^{12}=1$ that are studied in mirror symmetry, \cite{coxkatz}.

Dwork and Bombieri have conjectured conversely that all G-operators come from geometry. So one may ask: is the local system of solutions
$Sol(P)$ of a CY(n)-operator of the form $\C \otimes \mathbb{L}$, where $\mathbb{L}$ is a CY(n)-local system in the above sense? When can one
 achieve $\mathbb{L}=\H$? If $\mathbb{L}=\H$, can one find a family $f:X\lra \P^1$ with generic fibre a {\em Calabi-Yau} $n-1$-fold?\\
We refer to \cite{EvS} for a conjectural approach via mirror symmetry for CY(4)-operators.\\

Now suppose the whole situation is defined over $\Z$ and consider the reduction of $X \lra \P^1$ modulo some prime number $p$.
The object ${\cal L} \subset R^{n-1}f_*(\Q_l)$ ($l \neq p$ now defines an $l$-adic sheaf on
$\P^1$, lisse (that is, smooth) in some subset $S$. In particular, for each point $s :Spec(k) \lra S$ one has an action of
$Gal(\overline{k}/k)$ on the stalk ${\cal L}_s$. Hence one obtains a Frobenius element $Frob_s \in Aut({\cal L}_s)$
and 
\[P_{s}(T):=\det(1-T.Frob_s) \in \Z[T]\]
determines a factor of the zeta function of the reduction $X_s\mod p$.\\
 
To get a computational handle on these Frobenius polynomials it turns out to be useful to  change to a de Rham-type
description of the cohomology, \cite{kedlaya}.  It was Dwork who realised early that there is a tight interaction between the Gauss-Manin 
connection and the Frobenius operator. This leads in general to a relation between periods and the zeta function and 
in 1958 he gave his famous  $p-$adic analytic formula for the Frobenius polynomial in 
terms of a solution of the Picard-Fuchs differential equation for the Legendre family of elliptic curves, which we will now review.\\
The affine part of the Legendre family is given by
\[X_s: y^2=x(x-1)(x-s),\]
where $s\not=0,1$.
Over $\mathbb{C}$, the relative de Rham cohomology $H^1_{dR}$ of the family is free of rank 2, and the Hodge filtration $\textnormal{Fil}^1H^1_{dR}$ is generated by the differential 
\[\omega:=\frac{dx}{y}.\]
Let $\nabla$ be the Gauss-Manin connection on $H^1_{dR}$. Then, $\omega$ satisfies the differential equation
\[s(s-1)\omega''+(1-2s)\omega'-\frac{1}{4}\omega=0,\]
where $\omega'=\nabla(\omega)$. Let $f_0$ be the unique solution in $\mathbb{C}[[s]]$ to the above differential equation satisfying $f_0(s)=1$. $f_0$ is then given by the hypergeometric series
\[f_0(s)=F\left(\frac{1}{2},\frac{1}{2},1,s\right)=\sum_{j=0}^{\infty}\left(\frac{\left(\frac{1}{2}\right)_j}{j!}\right)^2s^j.\]
 
Now let $s_0 \in \mathbb{F}_{p^a}$ such that $f_0^{(p-1)/2}(s_0)\not=0$, where $f_0^{(p-1)/2}$ is the truncation of $f_0$ up to degree $(p-1)/2$. Let $\hat{s}$ be the Teichm\"uller lifting of $s_0$ to $W(\mathbb{F}_{p^a})$. The formal power series 
\[h(s):=\frac{f_0(s)}{f_0(s^p)}\]
converges at $\hat{s}$ and can be evaluated there. If $\epsilon =(-1)^{(p-1)/2}$, the element
\[\pi:=\epsilon^af_0(\hat{s})f_0(\hat{s}^p)...f_0(\hat{s}^{p^{a-1}})\]
is a reciprocal zero of the Frobenius polynomial, and the zeta function of $X_{s_0}$ is given by
\[\zeta(X_{s_0},T)=\frac{(1-\pi T)(1-p^a/\pi T)}{(1-T)(1-p^aT)}.\]

Thus, Dwork found a way to derive a formula for the Frobenius polynomial, which does only depend (up to $\epsilon$) on the solution of the Picard-Fuchs differential equation. The geometrical origin of $\epsilon$ lies in the geometry of the
singular fibre  $X_0$, which has a node with tangent cone $x^2+y^2=0$, which splits over $\F_p$ precisely when $\epsilon=1$.\\

In this paper we will consider the following\\
 {\bf Question:} Given a CY(n)-operator $P$ of $f: X \lra \P^1$ defined over $\Z$, is there a way
to calculate the Frobenius polynomials $P_{s}(T)$?\\

We describe a method to solve this problem for $n=4$ (modulo ``$\epsilon$'') and illustrate the procedure on some non-trivial examples.

\section{Dworks Unit-root Crystals}

We give a short introduction to the theory of Hodge $F$-crystals, which provides a framework to formalise the interaction between the Gauss-Manin connection and the Frobenius operator. (see \cite{stienstra}, \cite{katz}, \cite{deligne1},\cite{yu}).\\ 

Let $k$ be a perfect field of characteristic $p>0$, and let $W(k)$ be the ring of Witt vectors of $k$. Let $A$ be the ring $W(k)[z][g(z)^{-1}]$, where $g$ is a polynomial in $z$ (which will be specified later according to the actual situation), and let $A_n$ be the ring $A/p^{n+1}A$. By $A_{\infty}:=\displaystyle\lim_{\leftarrow} A/p^{n+1}A$, we denote the $p-$adic completion of $A$.\\
Let $\sigma$ be the absolute Frobenius on $k$, given by $\sigma(x)=x^p$. Following \cite{yu}, we define
\begin{definition}\begin{enumerate}
\item An $F-$crystal over $W(k)$ is a free $W(k)-$module $H$ of finite rank with a $\sigma-$linear endomorphism
\[F:H\rightarrow H\] such that $F\otimes\mathbb{Q}_p:H\otimes\mathbb{Q}_p\rightarrow H\otimes\mathbb{Q}_p$ is an isomorphism. If $F$ iself is an isomorphism, we call $H$ a {\em unit-root} $F-$crystal.
\item A Hodge $F-$crystal over $W(k)$ is an $F-$crystal $H$ equipped with a filtration by free $W(k)-$ submodules
\[H=\textnormal{Fil}^0H\supset\textnormal{Fil}^1H\supset...\supset\textnormal{Fil}^{N-1}H\supset\textnormal{Fil}^{N}H=0\] (called the Hodge filtration on $H$) which satisfies $F(\textnormal{Fil}^iH)\subset p^iH$ for all $i$.
\end{enumerate}
\end{definition}

The Frobenius automorphism $\sigma$ on $k$ lifts canonically to an automorphism $\sigma$ on $W(k)$. \\ There are different lifts of the Frobenius $\sigma$ on $A_{\infty}$, which restrict to $\sigma$ on $W(k)$ and reduce to the $p-$th power map modulo $p$. Let $\phi:A_{\infty}\rightarrow A_{\infty}$ be such a lift of Frobenius. 
\begin{definition} An $F-$crystal over $A_{\infty}$ is  a finitely generated free $A_{\infty}-$module $H$ with an integrable and $p-$adically nilpotent connection
\[\nabla:H\rightarrow \Omega_{A_{\infty}/W(k)}\otimes_AH.\] such that for every lift $\phi:A_{\infty}\rightarrow A_{\infty}$ of Frobenius, there exists a homomorphism of $A_{\infty}-$ modules
\[F(\phi): \phi^*H\rightarrow H\] such that the square
\[\begin{xy}\xymatrix{H\ar[rr]^{\nabla}\ar[d]_{F(\phi)\phi^*}& &\Omega^1_{A_{\infty}/W(k)}\otimes H\ar[d]^{\phi\otimes F(\phi)\phi^*}\\ H\ar[rr]^{\nabla} & &\Omega^1_{A_{\infty}/W(k)}\otimes H}
\end{xy}\]
is commutative, and such that $F(\phi)\otimes \mathbb{Q}_p: \phi^*H\otimes\mathbb{Q}_p\rightarrow H\otimes\mathbb{Q}_p$ is an isomorphism.  If $F(\phi)$ itself is an isomorphism, we call $H$ a unit-root crystal.
\end{definition}
From now on, to simplify the notation, we set $F:=F(\phi)\phi^*$.

\begin{definition} A divisible Hodge $F-$crystal $H$ is an $F-$crystal $H$ equipped with a filtration by free $A_{\infty}-$submodules
\[H=\textnormal{Fil}^0H\supset\textnormal{Fil}^1H\supset...\supset\textnormal{Fil}^{N-1}H\supset\textnormal{Fil}^{N}H\] 
(called the Hodge filtration on $H$) which satisfies
\begin{enumerate}
\item $\nabla \textnormal{Fil}^iH\subset\Omega^1_{A_{\infty}/W(k)}\otimes_{A_{\infty}}\textnormal{Fil}^{i-1}H$
\item $F(\textnormal{Fil}^iH)\subset p^iH$.
\end{enumerate}
\end{definition}

\begin{proposition}\label{wedge2} Let $H$ be a divisible Hodge $F-$crystal where $H/\textnormal{Fil}^1H$ is free of rank one. Then $\wedge^2H$ is a divisible Hodge $F-$crystal, with homomorphism of $A_{\infty}-$modules 
\[\frac{1}{p}\wedge^2F:\wedge^2H\rightarrow \wedge^2H\]and with Hodge filtration  given by \[\textnormal{Fil}^{i-1}(\wedge^2H)=\sum_{k=0}^i\textnormal{Fil}^kH\wedge\textnormal{Fil}^{i-k}H\]
for $i \geq 1$.
\end{proposition}
{\em Proof:} \ \ 
Since $H/\textnormal{Fil}^1H$ is of rank one, $\textnormal{Fil}^0\wedge\textnormal{Fil}^0=\textnormal{Fil}^0\wedge \textnormal{Fil}^1$.\\
Let $a \in \textnormal{Fil}^kH$ and $b \in \textnormal{Fil}^{i-k}H$.Then, $a\wedge b \in \textnormal{Fil}^{i-1}(\wedge^2 H)$ and \[\frac{1}{p}\wedge^2F(a\wedge b)=\frac{1}{p}Fa\wedge Fb \in \frac{1}{p}p^kH\wedge p^{i-k}H=p^{i-1}\wedge^2 H.\]
For $i \geq 2$,
\[\nabla(a\wedge b)=\nabla(a)\wedge b+a\wedge \nabla (b) \in \Omega_{A_{\infty}/W(k)}\otimes_{A_{\infty}}\textnormal{Fil}^{i-2}(\wedge^2 H).\]
\vspace{2mm}

Let $k'$ be a perfect field extension of $k$ and let $e_0:A_0\rightarrow k'$ be a $k-$morphism. Let $e_0(z)=\alpha_0$,  let $\alpha$ be the Teichm\"uller lifting of $\alpha_0$ in $W(k')$ and let $e:A_{\infty}\rightarrow W(k')$ be the $W(k)-$morphism with $e(z)=\alpha$. By $H_{\alpha}$, we denote the Teichm\"uller representative $H_{\alpha}:=H\otimes_{(A_{\infty},e)}W(k')$ of the crystal $H$ at the point $e_0$, which is an $F-$crystal with corresponding map $F_{\alpha}:=e^*F$. If $H$ is a Hodge $F-$crystal, then so is $H_{\alpha}$.\\
On $W(k')[[z-\alpha]]$, we put the natural connection $\nabla$ and choose the lift of Frobenius given by $\phi(z)=z^p$.
\begin{theorem}\label{katz41}(\cite{yu}, Theorem 2.1 or \cite{katz}, Theorem 4.1)\\ 
Let $\bar{k}$ be the algebraic closure of $k$, and let $H$ be a divisible Hodge $F-$crystal over $A_{\infty}$.\\
 If $H/\textnormal{Fil}^1H$ is of rank one and if for every $k-$morphism $e_0:A_0\rightarrow \bar{k}$ with $e_0(z)=\alpha_0$ and $\alpha \in W(\bar{k})$ a Teichm\"uller lifting of $\alpha_0$, $H_{\alpha}$ contains a direct factor of rank one, transversal to $\textnormal{Fil}^1H_{\alpha}$, which is fixed by the map induced by $F$ on $H_{\alpha}$, then there exists a unique unit-root F-subcrystal $U$ of $H$ such that $H=U\oplus \textnormal{Fil}^1H$ as $A_{\infty}$ modules.\\
 Suppose that over $A_{\infty}$, $U$ is locally generated by $u$. Write $F(u)=r(z)u$ for $r(z) \in A_{\infty}^*$. Then we have
\begin{enumerate}
\item Let $e_0:A_0\rightarrow k'$ be a $k-$morphism to a perfect field extension $k'$ of $k$ with $e_0(z)=\alpha_0$ where $u$ is defined. Let $\alpha$ be the Teichm\"uller lifting of $\alpha_0$. Then there exists an $f_0 \in W(k')[[z-\alpha]]$ such that $v:=f_0\cdot u \in W(k')[[z-\alpha]]\otimes_{A_{\infty}}H$ is horizontal with regard to $\nabla$ and the quotient $f_0/f_0^{\phi}$ is in fact the expansion of an element in $A_{\infty}$. 
\item There exists $c \in W(\bar{k})$ such that $c\cdot v \in W(\bar{k})\otimes_{W(k)}H$ is fixed by $F$ and $r(z)=(cf_0)/(cf_0)^{\phi}$.
\end{enumerate}
\end{theorem} 
 
The fact that $f_0/f_0^{\phi} \in A_{\infty}$ although $f_0 \in W(k)[[z-\alpha]]$ means that $f_0/f_0^{\phi}$ is a local expression of a ``global'' function. Although $f_0$ itself does only converge in a neighbourhood of $\alpha$, the global function expressed by the ratio $f_0/f_0^{\phi}$ converges at any Teichm\"uller point in Spec$(A_{\infty})$.

\subsection{CY(4)-operators and the corresponding crystals}
Now let ${\cal P}$ be a CY(4)-operator. We assume that ${\cal P}$ is the Picard-Fuchs operator  on a rank 4 submodule $H \subset H^3_{dR}(X/S_{\infty})$ for some family $f:X\rightarrow S_{\infty}$ of smooth CY-threefolds.\\
Let $k$ be the finite field with $p^r$ elements.
 From now on, we have $S_{\infty} =$ Spec$(A_{\infty})$, where $A=W(k)[z][(zs(z)g(z))^{-1}]$ for some polynomials $g(z)$ and $s(z)$. We assume that over the roots of $s(z)$, the family becomes singular. We will specify the polynomial $g(z)$ later (see section \ref{polyg}); it will be chosen in a way such that over each Teichm\"uller point $\alpha \in S_{\infty}$, the Frobenius polynomial on $H_{\alpha}$ is of the form

\[P:=1+aT+bpT^2+ap^3T^3+p^6T^4.\]
with four different reciprocal roots 
\[r_1,pr_2,p^2/r_2,p^3/r_1,\]
where $r_1$ and $r_2$ are $p-$adic units. Hence, giving a formula for the polynomial $P$ is equivalent to giving formulas for the $p-$adic units $r_1$ and $r_2$.\\
In general, if $f:V\rightarrow V$ is a homomorphism of vector spaces, then the eigenvalues of $\wedge^2f:\wedge^2V\rightarrow\wedge^2V$ are given by products $ab$, where $a$ and $b$ are eigenvalues of $f$ corresponding to linearly independent eigenvectors. \\
Let $\alpha_0 \in S_0$, and let $\alpha \in S_{\infty}$ be the Teichm\"uller lifting of $\alpha_0$.
 By Proposition \ref{wedge2}, the Frobenius automorphism on each fiber $\wedge^2H_{\alpha}$ of the crystal $\wedge^2H$ is given by $\frac{1}{p}\wedge^2F_{\alpha}$, where $F_{\alpha}$ is the Frobenius on $H_{\alpha}\subset H_{dR}^3(X_{\alpha})$. The eigenvalues of the relative Frobenius $(\wedge^2F_{\alpha})^r$ on the fibres  $\wedge^2H_{\alpha}$ are of the form $a_{\alpha}b_{\alpha}/p$, where $a_{\alpha}$ and $b_{\alpha}$ are eigenvalues of the relative Frobenius $F_{\alpha}^r$ on the corresponding fibre $H_{\alpha}$. Thus, if $r_1$ is the unit root  on a fibre  $H_{\alpha}$, and $\hat{r}_1$ is the unit root on the corresponding fibre  $\wedge^2H_{\alpha}$, then the roots of the Frobenius polynomial $\det(1-TF_{\alpha}^r)$ on $H_{\alpha}$ are given by
\begin{equation}r_1,p\hat{r}_1/r_1,p^2r_1/\hat{r}_1,p^3/r_1.\label{roots}\end{equation}  
We will give $p-$adic analytic formulas for the unit roots $r_1$ and $\hat{r}_1$.\\

\subsection{Horizontal sections for CY differential operators of order 4 and 5}\label{hors}
Let ${\cal P}$ be a CY(4)-operator. The differential equation 
 ${\cal P}y=0$ can be written in the form
\[y^{(4)}+a_3y^{(3)}+a_2y^{(2)}+a_1y^{(1)} +a_0y=0,\]
 where the coefficients $a_i$ satisfy the following relation:
\begin{equation}a_1=\frac{1}{2}a_2a_3-\frac{1}{8}a_3^3+a_2'-\frac{3}{4}a_3a_3'-\frac{1}{2}a_3''.\label{CY2}\end{equation}
\begin{proposition}(see \cite{yu})\label{su4}
Let ${\cal P}$ be a CY(4) differential operator and let $(H,\nabla)$ be a $\mathbb{Q}(z)/\mathbb{Q}$ differential module. Let $\omega \in H$ such that
\[\nabla^{4}\omega+a_3\nabla^3\omega+a_2\nabla^2\omega+a_1\nabla\omega +a_0\omega=0\] and let $f_0 \in \mathbb{Q}[[z]]$ be a formal solution to the differential equation ${\cal P}y=0$. If  $Y:=\exp\left(1/2\int a_3\right) \in \mathbb{Q}[[z]]$, then the following element $u_4 \in H\otimes_{\mathbb{Q}[z]}\mathbb{Q}[[z]]$ is horizontal with regard to $\nabla$ :
 \begin{eqnarray}u_4 & = & Y[f_0\nabla^3(\omega)-f_0'\nabla^2(\omega)+f_0''\nabla(\omega)-f_0'''\omega]+(Y a_3-Y')[f_0\nabla^2(\omega)-f_0''\omega]\nonumber \\ &+ &(Y a_2-(Y a_3)'+Y'')[f_0\nabla(\omega)-f_0'\omega]\label{u4}.\end{eqnarray} \end{proposition}
{\em Proof:} See \cite{yu}. The proof is by direct computation, using (\ref{CY2}).

Now let ${\cal Q}$ be a CY(5)-operator. The differential equation ${\cal Q}y=0$ can be written in the form
\[y^{(5)}+b_4y^{(4)}+b_3y^{(3)}+b_2y^{(2)}+b_1y^{(1)}+b_0y=0.\] 

\begin{proposition} The operator ${\cal Q}$ satisfies the second condition for CY(5) of the introduction, if and only if
the coefficients $b_i(z)$ satisfy the relations
\begin{equation}b_2=\frac{3}{5}b_3b_4-\frac{4}{25}b_4^3+\frac{3}{2}b_3'-\frac{6}{5}b_4b_4'-b_4''\label{b2}\end{equation}
and
\begin{eqnarray}
b_0 & = &\frac{1}{2}b_1'-\frac{2}{125}b_3b_4^3+\frac{1}{5}b_1b_4-\frac{1}{10}b_3b_4''+\frac{2}{5}b_4'''b_4+\frac{4}{5}b_4''b_4'+\frac{16}{125}b_4'b_4^3\\ \nonumber &+&\frac{12}{25}(b_4')^2b_4 - \frac{3}{10}b_3''b_4 + \frac{8}{25}b_4^2b_4''-\frac{3}{10}b_3'b_4'-\frac{3}{25}b_4^2b_3'-\frac{1}{4}b_3'''+\frac{16}{3125}b_4^5\\ \nonumber &+ &\frac{1}{5}b_4''''-\frac{3}{25}b_3b_4'b_4.\label{b0}
\end{eqnarray} 
\end{proposition}
{\em Proof:} By direct calculation, for details we refer to \cite{bogner}. 

\begin{proposition}\label{su5}
Let ${\cal Q}$ be a CY(5) differential operator and let $(H,\nabla)$ be a $\mathbb{Q}(z)/\mathbb{Q}$ differential module. Let $\eta \in H$ such that
\[\nabla^5\eta+ b_4\nabla^{4}\eta+b_3\nabla^3\eta+b_2\nabla^2\eta+b_1\nabla\eta +b_0\eta=0\] and let $F_0 \in \mathbb{Q}[[z]]$ be a formal solution to the differential equation ${\cal Q}y=0$. If $Y:=\exp\left(2/5\int b_4\right) \in \mathbb{Q}[[z]]$, then the following element $u_5 \in H\otimes_{\mathbb{Q}[z]}\mathbb{Q}[[z]]$ is horizontal with regard to $\nabla$:
\begin{eqnarray}
u_5 & = &Y [F_0\nabla^4(\eta)-F_0'\nabla^3(\eta)+F_0''\nabla^2(\eta)-F_0'''\nabla(\eta)+F_0''''\eta]\nonumber\\ &+& (Y b_4-Y')[F_0\nabla^3(\eta)-\frac{1}{3}F_0'\nabla^2(\eta)-\frac{1}{3}F_0''\nabla(\eta)+F_0'''\eta]\nonumber\\
 & +&  (Y b_3-(Y b_4)'+Y'')[F_0\nabla^2(\eta)+F_0''\eta]+(\frac{4}{3}((Y b_4)'-Y'')-\alpha b_3)F_0'\nabla(\eta)\nonumber\\
& + & (\frac{1}{2}((Y b_3)'-\frac{4}{3}((Y b_4)''-Y''')))[F_0'\eta+F_0\nabla(\eta)]\nonumber\\ &+ &(Y b_1-\frac{1}{2}((Y b_3)'-\frac{4}{3}((Y b_4)'-Y'''')))F_0\eta\label{u5}, 
\end{eqnarray}
 
\end{proposition}
{\em Proof:} Applying the identities (\ref{b2}) and (\ref{b0}), one directly verifies that $u_5$ satisfies $\nabla(u_5)=0$.\\

\subsection{Dwork's congruences}
Let $p$ be a prime number. We say that a sequence $(c_n)_{n \in \N}$ satisfies the {\em Dwork-congruences for $p$}, if the associated sequence
$C(n):=c(n)/c(\lfloor \frac{n}{p} \rfloor )\in \mathbb{Z}_p$ satisfies
\[ C(n) \equiv C(n+mp^s) \mod p^s\] 
for all $n,s \in N$ and $m=\{0,1,\ldots,p-1\}$ and if $c(0)=1$. We say that the Dwork-congruences hold for a CY(n) differential operator ${\cal P}$ if the Dwork-congruences hold
 for the sequence  $(c_n)_{n \in \mathbb{N}}$ of coefficients of the holomorphic solution  \[f_0(z)=\sum_{n=0}^{\infty} c_n z^n\] to the differential equation ${\cal P}y=0$ around $z=0$.
Dwork shows (see \cite{dwork}, Corollary 1. and 2.) that hypergeometric type numbers 
satisfy these Dwork congruences for all $p$.\\

\begin{theorem}(see \cite{dwork}, Lemma 3.4.)
Let $y(z)=\sum_nc_nz^n$ such that $(c_n)$ satisfies the Dwork congruences. Let $D:=\{x \in \mathbb{Z}_p, |y^{(p-1)}(x)|=1\}$. Then, for all $x \in D$, 
\[\frac{y(z)}{y(z^p)}|_{z=x} \equiv \frac{y^{(p^s-1)}(x)}{y^{(p^{s-1}-1)}(x^p)} \mod p^s.\]
\end{theorem}

This leads to an efficient evaluation of the left hand side at Teichm\"uller points. 
(Here
$y^{(p^s-1)}(z)$ is the polynomial obtained from $y(z)$ by truncation at $z^{p^s}$.) This crucial fact
was used in all of our computations.

\subsection{A formula for the roots of the Frobenius polynomial}\label{polyg}
Let ${\cal P}:={\cal P}(\theta,z)$ be a CY(4)-operator, where $\theta$ denotes the logarithmic derivative $z\partial/\partial z$. \\
As before, we assume that ${\cal P}$ is the Picard- Fuchs operator  on a rank four submodule $H\subset H_{dR}^3(X/S_{\infty})$ for a family $f:X\rightarrow S_{\infty}$ of smooth CY- threefolds.\\ 
The rank 6$={4 \choose 2}$ $A_{\infty}-$ module $\wedge^2H$ is a direct sum of an $A_{\infty}-$ module $G$ of rank 5 and a rank 1 module. The rank 1 module is generated by a section that corresponds to the pairing $<-,->$ and is horizontal with respect to $\nabla$.\\ 
We construct a 5th order differential operator ${\cal Q}$ on the submodule $G$ by choosing ${\cal Q}$ to  be the differential operator of minimal order such that for any two linearly independent solutions $y_1(z),y_2(z)$ of the differential equation ${\cal P}y=0$, 
\[w:=z\left|\begin{array}{cc}y_1 & y_2 \\ y_1' & y_2'\end{array}\right|\]
is a solution of ${\cal Q}w=0$ .

\begin{proposition}
The operator ${\cal Q}$ satisfies the first and the second condition of CY(5).
\end{proposition}
{\em Proof:}
The statement that ${\cal Q}$ satisfies the first condition of CY(5) is the content of \cite{almkvistzudilin}, Proposition 4. 
A direct computation shows that since ${\cal P}$ is a $CY(4)-$operator, the coefficients of ${\cal Q}$ satisfy the
equations (\ref{b0}) and (\ref{b2}), so the second condition of CY(5) holds.\\

In all examples it was found that the operator ${\cal Q}$ also has an integral power series solution, and thus satisfies the third condition of CY(5).
For the moment, however, we are unable to prove this is general so we  
\begin{conjecture}\label{CY5Q}
The differential operator ${\cal Q}$, constructed from a CY(4)-operator ${\cal P}$ as above, satisfies the third condition of CY(5).
\end{conjecture}
So if Conjecture \ref{CY5Q} holds true, the differential operator ${\cal Q}$ is a CY(5)-operator.\\

${\cal Q}$ can be expressed in terms of $\wedge^2{\cal P}(\theta,z)$ as
\[{\cal Q}(\theta,z)=\wedge^2{\cal P}(\theta-1,z).\]
For the differential operators ${\cal P}$ and ${\cal Q}$, we use the same notation with coefficients $a_i$ and $b_i$ as in section \ref{hors}.
\begin{proposition}\label{omegaeta}
Let ${\cal Q}$ be the CY(5)-operator constructed above, and let $\omega \in H$ such that 
\[\nabla^{4}\omega+a_3\nabla^3\omega+a_2\nabla^2\omega+a_1\nabla\omega +a_0\omega=0.\]
Then, the element $\eta:=z\omega\wedge\nabla\omega \in G$ satisfies
\[\nabla^5\eta+ b_4\nabla^{4}\eta+b_3\nabla^3\eta+b_2\nabla^2\eta+b_1\nabla\eta +b_0\eta=0.\]
\end{proposition}
{\em Proof:} The proposition follows by a straightforward calculation, applying the relations between the coefficients $a_i$ of the CY(4)-operator ${\cal P}$ and the coefficients $b_i$ of the CY(5)-operator ${\cal Q}$ listed in \cite{almkvist}.\\

It still remains to point out how to choose the polynomial $g(z)$ in the definition of the ring $W(k)[z][(zs(z)g(z))^{-1}]$ to obtain divisible Hodge $F-$crystals $H\subset H_{dR}^3(X/S_{\infty})$ and $G \subset \wedge^2H$ which satisfy the conditions of Theorem \ref{katz41}.\\

The following conjecture was crucial for the choice of the polynomial $g(z)$:

\begin{conjecture}\label{keyconj}\ \ \
\begin{enumerate}
\item Let $f_0$ be the solution of the differential equation ${\cal P}y=0$ around $z=0$ with $f_0(0)=1$. If the coefficients $c_n$ in the expansion 
\[f_0(z)=\sum_{n=0}^{\infty}c_nz^n\]
satisfy the Dwork congruences, then $H$ satisfies the  conditions of Theorem (\ref{katz41}) if the polynomial $g(z)$ in the definition of $A_{\infty}$ is chosen as $g(z):=f^{(p-1)}_0(z)$.
\item Let $F_0(z)$  be the solution of the differential equation ${\cal Q}y=0$ around $z=0$ with $F_0(0)=1$. If the coefficients $d_n$ in the expansion 
\[F_0(z)=\sum_{n=0}^{\infty}d_nz^n\]
satisfy the Dwork congruences, then the sub-$F$- crystal $G \subset \wedge^2H$ satisfies the  conditions of Theorem (\ref{katz41}) if the polynomial $g(z)$ in the definition of $A_{\infty}$ is chosen as $g(z):=F^{(p-1)}_0(z)$.
\end{enumerate} 
\end{conjecture}

According to the conjecture, it seems to be the right thing to choose $g(z)=f_0^{(p-1)}(z)F_0^{(p-1)}(z)$. So from now on, we fix the ring $A_{\infty}$ by 
\[A:=W(k)[z][(zs(z)f_0^{(p-1)}(z)F_0^{(p-1)}(z))^{-1}].\]
This choice was confirmed by our numerous computations; for each parameter value $z=\alpha$ with $f_0^{(p-1)}(\alpha_0) \not= 0 \mod p$ and $F_0^{(p-1)}(\alpha_0) \not= 0 \mod p$, in the examples we considered, we were able to compute the Frobenius polynomial explicitly.\\

For each pair of CY(4) and CY(5) operators we treat in this paper, the functions
\[Y_4=\exp\left(1/2 \int a_3\right)\; \textnormal{and} \; Y_5=\exp\left(2/5\int b_4\right)\] satisfy $Y_4\in \mathbb{Q}(z)$ and $Y_5 \in \mathbb{Q}(z)$.
Thus, in each of the examples we considered, the following proposition holds:
\begin{proposition}
Let \[r(z)=\frac{f_0(z)}{f_0(z^p)} \; \textnormal{and} \; \hat{r}(z)=\frac{F_0(z)}{F_0(z^p)}.\]
Assuming that Conjecture \ref{keyconj} holds,  if $|\alpha s(\alpha)f_0^{(p-1)}(\alpha)F_0^{(p-1)}(\alpha)|=1$, there exist constants $
\epsilon_4$ and $\epsilon_5 \in W(\bar{k})$ such that the $p-$ adic units $r_1(\alpha)$ and $\hat{r}_1(\alpha)$ determining the Frobenius polynomial on $H_{\alpha} \subset H_{dR}^3(X_{\alpha})$ are given by  
\[r_1(\alpha)=(\epsilon_4^{(1-\sigma)})^{1+...+\sigma^{r-1}}r(\alpha)r(\alpha^p)...r(\alpha^{p^{r-1}})\]
and
\[\hat{r_1}(\alpha)=(\epsilon_5^{(1-\sigma)})^{1+\sigma+...+\sigma^{r-1}}\hat{r}(\alpha)\hat{r}(\alpha^p)....\hat{r}(\alpha^{p^{r-1}}).\]
If we assume furthermore that
\begin{equation*}(\epsilon_4^{1-\sigma})^{1+\sigma+...+\sigma^{{r-1}}}=(\epsilon_5^{1-\sigma})^{1+\sigma+...+\sigma^{r-1}}=1,\label{constants}\end{equation*} the $p-$adic units are given by
\[r_1(\alpha)=r(\alpha)r(\alpha^p)...r(\alpha^{p^{r-1}})\]
and
\[\hat{r_1}(\alpha)=\hat{r}(\alpha)\hat{r}(\alpha^p)....\hat{r}(\alpha^{p^{r-1}}).\]
\end{proposition} 
{\em Proof:}\ \ There exists an $\omega \in H$ such that the horizontal section w.r.t. $\nabla$ is given by formula (\ref{u4}), while on $G$, it is given by formula (\ref{u5}), where $\eta =z\omega\wedge\nabla\omega$ by Proposition \ref{omegaeta}. These sections $u_4$ and $u_5$ play the role of the section $v=f\cdot u$ in Theorem \ref{katz41}. Hence, the section $u$ in the theorem is given by $(f_0Y_4)^{-1}u_4$ and $(F_0Y_5)^{-1}u_5$ respectively, where $Y_4=\exp\left(1/2 \int a_3\right)\in \mathbb{Q}(z)$ and $Y_5=\exp\left(2/5\int b_4\right)\in \mathbb{Q}(z)$.  Since 
\[\left(\frac{Y_4(z)}{Y_4(z^p)}|_{z=\alpha}\right)^{1+\sigma+...+\sigma^{r-1}}=1 \; \textnormal{and} \;\left(\frac{Y_5(z)}{Y_5(z^p)}|_{z=\alpha}\right)^{1+\sigma+...+\sigma^{r-1}}=1,\]
by Theorem \ref{katz41} there exist constants $\epsilon_4$ and $\epsilon_5 \in W(\bar{k})$ (where $\bar{k}$ denotes the algebraic closure of $\mathbb{F}_p$) such that
\[r_1(\alpha)=(\epsilon_4^{(1-\sigma)})^{1+...+\sigma^{r-1}}r(\alpha)r(\alpha^p)...r(\alpha^{p^{r-1}})\]
and
\[\hat{r_1}(\alpha)=(\epsilon_5^{(1-\sigma)})^{1+\sigma+...+\sigma^{r-1}}\hat{r}(\alpha)\hat{r}(\alpha^p)....\hat{r}(\alpha^{p^{r-1}}).\]
Now we assume that the constants satisfy
\begin{equation}(\epsilon_4^{1-\sigma})^{1+\sigma+...+\sigma^{{r-1}}}=(\epsilon_5^{1-\sigma})^{1+\sigma+...+\sigma^{r-1}}=1.\label{constants}\end{equation}
Then, the proposition follows.

\section{Some special Picard-Fuchs equations}

We will apply the method explained in the previous section to compute Frobenius polynomials for some special
fourth order operators. These operators belong to the list \cite{AESZ}. A typical example is operator $45$ from that
list:
\[{\theta}^{4}-4\,x \left( 2\,\theta+1 \right) ^{2} \left( 7\,{\theta}^{
2}+7\,\theta+2 \right) -128\,{x}^{2} \left( 2\,\theta+1 \right) ^{2}
 \left( 2\,\theta+3 \right) ^{2}
\]
This operator is a so-called {\em Hadamard product} of two second order operators.

\subsection{Hadamard Products}

The {\em Hadamard product} of two power series $f(x):=\sum_n a_n x^n$ and $g(x)=\sum_n b_n x^n$ is the 
power-series defined by the coefficient-wise product:
\[f*g(x):=\sum_n a_n b_n x^n\]
It is a classical theorem, due to Hurwitz, that if $f$ and $g$ satisfy linear differential equations ${\cal P}$
and ${\cal Q}$ resp., then $f*g$ satisfies a linear differential equation ${\cal P}*{\cal Q}$. 
Only in very special cases, the Hadamard product of two CY-operators will again be CY, but it is a general fact that
if $f$ and $g$ satisfy differential equations of {\em geometrical origin}, then so does $f*g$. 
For a proof, we refer to \cite{andre}. Here we sketch the idea. The multiplication map
\[ m: \C^* \times \C^* \lra \C^*, (s,t) \mapsto s.t\] 
can be compactified to a map
\[\mu: \widetilde{\P^1 \times \P^1} \lra \P^1\]
by blowing-up the two points $(0,\infty)$ and $(\infty,0)$ of $\P^1 \times \P^1$. Given two families 
$X \lra \P^1$ and $Y \lra \P^1$ over $\P^1$, we define
a new family $X*Y \lra \P^1$,  as follows.
The cartesian product $X \times Y$ maps to $\P^1 \times \P^1$ and can be pulled back to $X*Y$ over 
$\widetilde{\P^1 \times \P^1}$. Via the map $\mu$ we obtain a family over $\P^1$. If
$n$ resp. $m$ is the fibre dimension of $X\lra \P^1$ resp. $Y \lra \P^1$, then $X*Y \lra \P^1$ has
fibre dimension $n+m+1$. The local system $H^{n+m+1}$ of $X*Y \lra \P^1$ contains the {\em convolution} of
the local systems of $X \lra \P^1$ and $Y \lra \P^1$. Note that the
critical points of $X*Y \lra \P^1$ are, apart from $0$ and $\infty$, the products of the critical values of the factors.
In down-to-earth terms, if $X \lra \P^1$ and $Y \lra \P^1$ are defined by say Laurent polynomials $F(x)$
and $G(y)$ resp., then the fibre of $X*Y \lra \P^1$ over $u$ is defined by the equations
\[F(x)=s, G(y)=t, s.t=u \] 
If the period functions for $X \lra \P^1$ and $Y \lra \P^1$ are represented as
\[f(s) =\int_\gamma Res(\frac{\omega}{F(x)-s})=\sum_n a_n s^n\]
\[g(t) =\int_\delta Res(\frac{\eta}{G(y)-t})=\sum_m b_m t^m\]
then 
\begin{eqnarray*}\int_{\gamma \times \delta \times S^1}\frac{\omega \wedge \eta \wedge ds \wedge dt}{(F(x)-s)(G(y)-t)(st-u)}&=&
\int_{S^1}\sum a_n s^n b_m t^m \frac{du}{u}\\ &= &\sum a_n b_n u^n=f(u)*g(u) \end{eqnarray*}
is a period of $X*Y \lra \P^1$.

For example, if we apply this construction to the rational elliptic surfaces $X=Y$ with singular fibres of Kodaira type
$I_9$ over $0$ and $I_1$ over $\infty$ and two further fibres of type $I_1$, we obtain a family $X * Y \lra \P^1$, with 
generic fibre a Calabi-Yau 3-fold with $h^{12}=1$ and $\chi=164$.

\subsection{Some special $CY(2)$-operators}

We will  use Hadamard-products of some very special CY(2)-operators appearing in $\cite{almkvistzudilin}$ from which
we also take the names. These operators all are associated to {\em extremal rational elliptic surfaces} $X \lra \P^1$ 
with non-constant  j-function. Such a surface has three or four singular fibres, \cite{mirandapersson}.
The six cases with three singular fibres fall into four isogeny-classes
and each of these gives rise to a Picard-Fuchs operator of hypergeometric 
type (named A,B,C,D) and one obtained by performing a M\"obius transformation interchanging  $\infty$ with the singular point $\neq 0$ (named e,h,i,j).

\[
\begin{array}{|c|l|l|}
\hline
\textup{Name}&\textup{Operator}&a_n\\
\hline
A&\theta^2-4x(2\theta+1)^2&\frac{(2n)!^2}{n!^4}\\
\hline
B&\theta^2-3x(3\theta+1)(3\theta+2)&\frac{(3n)!}{n!^3}\\
\hline
C&\theta^2-4x(4\theta+1)(4\theta+3)&\frac{(4n)!}{(2n)!n!^2}\\
\hline
D&\theta^2-12x(6\theta+1)(6\theta+5)&\frac{(6n)!}{(3n)!(2n)!n!}\\
\hline
\end{array}
\]

\[
\begin{array}{|c|l|l|}
\hline
\textup{Name}&\textup{Operator}&a_n\\
\hline
e&\theta^2-x(32\theta^2+32\theta+12)+256x^2(\theta+1)^2&16^n\sum_{k}(-1)^k \binom{-1/2}{k}\binom{-1/2}{n-k}^2\\
\hline
h&\theta^2-x(54\theta^2+54\theta+21)+729x^2(\theta+1)^2&27^{n}\sum_{k}(-1)^{k}\binom{-2/3}{k}\binom{-1/3}{n-k}^{2}\\
\hline
i&\theta^2-x(128\theta^2+128\theta+52)+4096x^2(\theta+1)^2&64^{n}\sum_{k}(-1)^{k}\binom{-3/4}{k}\binom{-1/4}{n-k}^{2}\\
\hline
j& \theta^2-x(864\theta^2+864\theta+372)+18664x^2(\theta+1)^2&432^{n}\sum_{k}(-1)^{k}\binom{-5/6}{k}\binom{-1/6}{n-k}^{2}\\
\hline
\end{array}
\]

The six cases with four singular fibres are the Beauville surfaces (\cite{beauville}) and also
form four isogeny classes and lead to the six Zagier-operators, called (a,b,c,d,f,g). 

These are also of the form 
\[\theta^2-x(a\theta^2+a\theta+b)-cx^2(\theta+1)^2\]
but now the discriminant $1-ax-cx^2$ is not a square, so the operator has 
four singular points.
\[
\begin{array}{|c|l|l|}
\hline
\textup{Name}&\textup{Operator}&a_n\\
\hline
a& \theta^2-x(7\theta^2+7\theta+2)-8x^2(\theta+1)^2&\sum_{k}\binom{n}{k}^{3} \\
c& \theta^2-x(10\theta^2+10\theta+3)+9x^2(\theta+1)^2 & \sum_{k}\binom{n}{k}^{2}\binom{2k}{k}\\
g&  \theta^2-x(17\theta^2+17\theta+6)+72x^2(\theta+1)^2&\sum_{i,j}8^{n-i}(-1)^{i}\binom{n}{i}\binom{i}{j}^{3}\\
\hline
d&   \theta^2-x(12\theta^2+12\theta+4)+32x^2(\theta+1)^2&\sum_{k}\binom{n}{k}\binom{2k}{k}\binom{2n-2k}{n-k}\\
\hline 
f& \theta^2-x(9\theta^2+9\theta+3)+27x^2(\theta+1)^2&\sum_{k}(-1)^{k}3^{n-3k}\binom{n}{3k}\frac{(3k)!}{k!^{3}}\\
\hline
b&  \theta^2-x(11\theta^2+11\theta+3)-x^2(\theta+1)^2&\sum_{k}\binom{n}{k}^{2}\binom{n+k}{n}\\ 
\hline
\end{array}
\]

The ten products $A*A$, etc. form 10 of the 14 hypergeometric families from \cite{AESZ}. 
The $16$ products $A*e$ etc. are not hypergeometric, but also have three singular fibres. The 24 operators
$A*a$ etc. have, apart from $0$ and $\infty$ two further singular fibres. The operators $a*a$ etc. have
four singular fibres apart from $0$ and $\infty$. \\

{\bf Observations:}\\ 
1) The Dwork-congruences hold for the operators $a,b,\ldots,j$. For the Apery-sequence (case b) this
was also conjectured in \cite{yu}. (It follows  from \cite{dwork} that $A,B,C,D$ satisfy the Dwork-congruences).
It follows that the Dwork-congruences hold for all fourth order Hadamard products within this group.

2) For the hypergeometric cases $A*A$ , etc, and the cases $A*a$, etc. the Dwork-congruences also hold for
the associated fifth order operator, although even for the simplest examples like the quintic threefold, this is not at all obvious. In the case of the quintic, the holomorphic solution around $z=0$ to the fifth order differential equation is given by the formula $F_0(z)=\sum_{n=0}^{\infty}A_nz^n$, where 
\[A_n:=\sum_{k=0}^n\frac{(5k)!}{k!^5}\frac{5(n-k)!}{(n-k)!^5}(1+k(-5H_k+5H_{n-k}+5H_{5k}-5H_{5(n-k)}))\]
and $H_k$ is the harmonic number $H_k=\sum_{j=1}^k\frac{1}{j}$. Thus, by the formula it is not even obvious that the coefficients $A_n$ are integers.

3) In fact, the Dwork-congruences hold for {\em almost all} fourth order operators from the list $\cite{AESZ}$.
It is an interesting problem to try to prove these experimental facts. On the other hand, it is clear that they 
cannot hold in general for differential operators of geometrical origin: if we multiply  $f_0$ with a rational 
function of $x$ we obtain a (much more complicated) CY-operator for which the congruences in general will not 
hold.

\subsection{Computations}
In the hypergeometric cases we reproduced results obtained in \cite{villegas}. In the appendix of \cite{ArxSvS}, the results of our 
calculations on the $24$ operators which are Hadamard products like $A*a$ etc. are collected.
We computed coefficients $(a,b)$ of the Frobenius polynomial 
\[P(T)=1+aT+bpT^2+ap^3T^3+p^6T^4 \]
for all primes $p$ between $3$ and $17$ and for all possible values of $z \in \mathbb{F}_p^*$. 
In our computations, we assumed that Conjecture \ref{keyconj} holds true and took the constants (\ref{constants}) appearing in the formula for the unit root to be one.
To generate the tables of coefficients in \cite{ArxSvS}, we used the programming language MAGMA. We computed with an overall $p-$adic accuracy of 500 digits. This was necessary, since in the computation of the power series solutions to the differential equations ${\cal P}y=0$ and ${\cal Q}y=0$, denominators divisible by large powers of $p$ occured during the calculations (although the solutions themselves have integral coefficients). The occurance of large denominators reduces the $p-$adic accuracy in MAGMA, and thus we had to compute with such a high overall accuracy to obtain correct results in the end. For the unit roots themselves, we computed the ratio
\[\frac{f_0(z)^{(p^3-1)}}{f_0(z^p)^{(p^2-1)}}|_{z=\alpha} \mod p^3\] with $p-$adic accuracy modulo $p^3$. 
We checked our results for the tuples $(a,b)$ determined the absolute values of the complex roots of the Frobenius polynomial, which by the Weil conjectures should have absolute value $p^{-3/2}$. Needless to say, this was always fulfilled. \\

\subsection{Example}
In this section, we describe the computational steps we performed in MAGMA for one specific example.
We consider the operator $A*a$, which is nr. 45 from the list \cite{AESZ}. \\
We compute the Frobenius polynomial for $p=7$ and $\alpha_0=2 \in \mathbb{F}_7$ with $4$ digits of $7-$adic precision, i.e. modulo $7^4$. Since $2 \not=-\frac{1}{16}$ and $2 \not=\frac{1}{128}$ in $\mathbb{F}_7$, $\alpha_0$ is not a singular point of the differential equation.\\
 First of all, we computed the truncated power series solution $f_0^{(p^{s+1}-1)}(z)$ to the differential equation 
\[{\cal P}y=0,\]
and obtained
\[f_0^{(7^4-1)}(z)=1+8z+360z^2+22400z^3+1695400z^4+143011008z^5+...\]
Thus,  $f_0^{(7-1)}(\alpha_0)= 1 \in \mathbb{F}_7$ is nonzero.
Let $\alpha^{(4)}$ be the Teichm\"uller lifting of $\alpha_0$ with $7-$adic accuracy of $4$ digits. Evaluating $f_0$ in this point, we obtain
\[f_0^{(7^4-1)}(\alpha^{(4)})\equiv 1709 \mod 7^4\]
and \[f_0^{(7^3-1)}((\alpha^{(4)})^7)\equiv 1814 \mod 7^4.\]
Thus, the unit root of the Frobenius polynomial is
\[r^4:=\frac{f_0^{(7^4-1)}(\alpha^{(4)})}{f_0^{(7^3-1)}((\alpha^{(4)})^7)} \equiv 582 \mod 7^4.\]  
To compute the second root of the Frobenius polynomial, we compute the truncated power series solution $F_0^{(7^4-1)}(z)$ of the fifth order differential equation
\[{\cal Q}y=0,\]
where ${\cal Q}$ is the second exterior power of the differential operator ${\cal P}$, given by
\begin{eqnarray*}{\cal Q} & = & \theta^5-z(44+260\theta+628\theta^2+792\theta^3+560\theta^4+224\theta^5)\\
&+&z^2(-6512+400\theta+44160\theta^2+71040\theta^3+42240\theta^4+8448\theta^5)\\
&+&z^3(4177920+13180928\theta+16588800\theta^2+10567680\theta^3+3440640\theta^4\\ &+&458752\theta^5)\\
&+&z^4(100663296+285212672\theta+310378496\theta^2+163577856\theta^3+41943040\theta^4\\ &+&4194304\theta^5). 
\end{eqnarray*}
The solution is given by
\[F_0^{(7^4-1)}=1+44z+3652z^2+337712z^3+33909700z^4+3567877424z^5 +...,\]
$F_0^{(7-1)}(\alpha_0)= 2 \in \mathbb{F}_7$ is nonzero
and we compute 
\[F_0^{(7^4-1)}(\alpha^{(4)})\equiv 51 \mod 7^4\]
and 
\[F_0^{(7^3-1)}((\alpha^{(4)})^7)\equiv 1387 \mod 7^4.\]
Thus, 
\[\hat{r}^4:=\frac{F_0^{(7^4-1)}(\alpha^{(4)})}{F_0^{(7^3-1)}((\alpha^{(4)})^7)} \equiv 1101 \mod 7^4.\]
Since the Frobenius polynomial (with $7-$adic accuracy 4) is given by
\[P(T)=(1-r^4T)(1-7\hat{r}^4/r^4T)(1-7^2r^4/\hat{r}^4T)(1-7^3/r^4T),\]
we finally obtain
\[P(T)=7^6T^4-7^3\cdot8T^3+7\cdot 2T^2-8T+1.\]
As expected, the complex roots of $P$ do have complex absolute value $7^{-3/2}$.\\

Exemplarily, we now list all values (a,b) we computed for the differential operator
$A*a$.
If there occurs a ``-'' in the table instead of the tuple $(a,b)$, then the correponding $z \in \mathbb{F}_p$ is either a zero of $f_0^{(p-1)}$ or $F_0^{(p-1)}$ or of both, where $f_0$ was the power series solution of the fourth order differential equation and $F_0$ was the solution of the fifth order equation. 
The appearance of $(a,b)'$ means that the polynomial is {\em reducible}. The
appearance of  $(a,b)^*$ means that the corresponding $z$ is a singular point of the
differential equation.\\
\newpage

\scriptsize
$p=3$ \ \ \ \ \ \ \ \ \ \ \ \ \ \ $p=5$:
\vspace{2mm}

$\begin{array}{|c|c|c|}\hline z&1&2 \\ \hline
&-&- \\ \hline
\end{array}$ \ \ \ \ 
$\begin{array}{|c|c|c|c|c|}\hline z&1&2&3&4 \\ \hline
&(6,-6)'&(28,38)*&-&(32,62)* \\ \hline
\end{array}$
\vspace{2mm}

\noindent $p=7$:
\vspace{2mm}

\begin{tabular}{|c|c|c|c|c|c|c|}\hline $z$&1&2&3&4&5&6 \\ \hline
&(2,-46)&(-8,2)&(32,-94)*&(80,290)*&(10,50)'&- \\ \hline
\end{tabular}
\vspace{2mm}

\noindent $p=11$:
\vspace{2mm}

\begin{tabular}{|c|c|c|c|c|c|c|c|c|}\hline $z$&1&2&3&4&5&6&7&8 \\ \hline
&(56,290)'&-&(-16,2)'&(6,26)&(16,98)&(12,114)'&(26,106)&- \\ \hline
\end{tabular}
\vspace{1mm}

\begin{tabular}{|c|c|c|}\hline $z$&9&10 \\ \hline
&(-8,2)&(-36,210)' \\ \hline
\end{tabular}
\vspace{2mm}

\noindent $p=13$:
\vspace{2mm}

\begin{tabular}{|c|c|c|c|c|c|c|c|c|}\hline $z$&1&2&3&4&5&6&7&8 \\ \hline
&(-8,270)'&(20,-106)&(-4,86)&(-204,646)*&(22,-30)&(-160,30)*&(-34,50)&(-16,302) \\ \hline
\end{tabular}
\vspace{1mm}

\begin{tabular}{|c|c|c|c|c|}\hline $z$&9&10&11&12 \\ \hline
&(58,146)&(18,34)&(84,406)&(56,206)' \\ \hline
\end{tabular}
\vspace{2mm}

\noindent $p=17$:
\vspace{2mm}

\begin{tabular}{|c|c|c|c|c|c|c|c|c|}\hline $z$&1&2&3&4&5&6&7&8 \\ \hline
&(256,-322)*&(256,-322)*&(-24,542)&(44,166)&(210,1218)&(24,-178)'&(-100,278)&(22,50) \\ \hline
\end{tabular}
\vspace{1mm}

\begin{tabular}{|c|c|c|c|c|c|c|c|c|}\hline $z$&9&10&11&12&13&14&15&16 \\ \hline
&(-4,70)&(52,470)&(-84,342)'&-&(22,-334)'&(18,258)&(184,974)&(-56,302)' \\ \hline
\end{tabular}
\vspace{2mm}
\normalsize

\subsection{Modular forms of weight four}

In some cases, the so chosen accuracy was too low, and we had to compute mod $p^4$. This happened in case the parameter $ \alpha_0 \in \mathbb{F}_p$ was  a critical point of the differential equation. But it is somewhat of a miracle that our
calculation made sense at the critical points at all. In order to understand what is supposed to happen at a singular point, 
recall that if the fibre $X_s$ of a family $X \lra \P^1$ over $s \in \P^1(\Q)$ aquires an ordinary double point, then 
the Frobenius polynomial should factor as
\[P(T)=(1-\chi(p)T)(1-p\chi(p)T)(1-a_pT+p^3T^2)\]
for some character $\chi$. The factor $(1-a_pT+p^3T^2)$ is the Frobenius polynomial on the two dimensional pure part 
of $H^3$. This part can be identified with the $H^3$ of a small resolution $\tilde{X_s}$, which then is a rigid Calabi-Yau $3$-fold.
According to the modularity conjecture for such Calabi-Yau $3$-folds, the coefficients $a_p$ are Fourier coefficients of a 
weight four modular form for some congruence subgroup $\Gamma_0(N)$, \cite{meyer}. 

This is exactly the phenomenon that occurs at the singular points of our differential equations.  For the hypergeometric 
cases we refind the results of \cite{villegas}. For $16$ of the $24$ operators $A*a$ etc, 
we have two rational critical values. In 31 of the cases we are able to identify the modular form.

We use the notation of modular forms as in \cite{meyer}: the notation $a/b$ means: the $b$-th Hecke
eigenform of level $a$. 'Twist of' means: the modular forms differ by character. We remark that
the critical points of the operators are reciprocal integers and the level of the corresponding 
modular form divides that integer. For the cases involving the operator $c$ 
one usually has equality and so the modular form for $D*c$ presumably has
level $3888$, which was outside the range of our table. Remark that all levels appearing
only involve primes $2$ and $3$.
$$
\begin{array}{|c||c|c|c||c|c|c|}
\hline
\textup{Case}&\textup{Point}&\textup{Form}&\textup{Twist of} &\textup{Point}&\textup{Form}&\textup{Twist of}\\
\hline
A*a&-1/16&8/1&-&1/128&64/5&8/1\\
B*a&-1/27&27/2&27/1&1/126&54/2&-\\
C*a&-1/64&32/3&32/2&1/512&256/3&-\\
D*a&-1/432&216/4&216/2&1/3456&1728/16&216/1\\
\hline
A*c&1/144&48/1&24/1&1/16&16/1&8/1\\
B*c&1/243&243/1&-&1/27&27/1&-\\
C*c&1/576&576/3&94/4&1/64&64/3&32/2\\
D*c&1/3888&&1944/5&1/432&432/9&216/2\\
\hline
A*d&1/128&64/4&32/1&1/64&32/2&-\\
B*d&1/216&9/1&-&1/108&108/4&108/2\\
C*d&1/512&256/1&-&1/256&128/4&128/1\\
D*d&1/3456&576/8&288/1&1/1728&864/3&864/1\\
\hline
A*g&1/144&24/1&-&1/128&64/1&8/1\\
B*g&1/243&243/2&243/1&1/216&54/4&54/2\\
C*g&1/576&288/10&96/4&1/512&256/4&256/3\\
D*g&1/3888&1944/6&1944/5&1/3456&1728/15&-\\
\hline
\end{array}
$$
\vspace{5mm}

The simplest modular forms appearing are the well-knonw $\eta$-products $8/1=\eta(q^2)^4\eta(q^4)^4$, $9/1=\eta(q^3)^8$.\\

{\bf Acknowledgement.}
We thank G. Almkvist and W. Zudilin for interest in the project. We thank
M. Bogner for help with calculations, S.Cynk and J. D. Yu  for useful discussions. Special thanks to V. Golyshev for his
idea to use ``polynomial functors'' to obtain information on the root of valuation one.\\
The work of K.S. was funded  by the SFB Transregio 45.

\newpage
\section{Appendix}
In this appendix we collect the results of our calculations on the $24$ operators Hadamard products $A*a$, etc.
We computed coefficients $(a,b)$ of the Frobenius Polynomial 
\[P(T)=1+aT+bpT^2+ap^3T^3+p^6T^4 \]
for all primes $p$ between $3$ and $17$ and for all possible values of $z \in \mathbb{F}_p^*$. If there occurs a ``-'' in the table instead of the tuple $(a,b)$, then the correponding $z \in \mathbb{F}_p$ is either a zero of $f_0^{(p-1)}$ or $F_0^{(p-1)}$ or of both, where $f_0$ was the power series solution of the fourth order differential equation and $F_0$ was the solution of the fifth order equation. 
The appearance of $(a,b)'$ means that the polynomial is {\em reducible}. The
appearance of  $(a,b)^*$ means that the corresponding $z$ is a singular point of the
differential equation.

\subsection{The Case $A*a$}

This is operator nr. 45 from the list \cite{AESZ}:
\[{\theta}^{4}-4\,x \left( 2\,\theta+1 \right) ^{2} \left( 7\,{\theta}^{
2}+7\,\theta+2 \right) -128\,{x}^{2} \left( 2\,\theta+1 \right) ^{2}
 \left( 2\,\theta+3 \right) ^{2}
\]
\scriptsize
$p=3$ \ \ \ \ \ \ \ \ \ \ \ \ \ \ $p=5$:
\vspace{2mm}

$\begin{array}{|c|c|c|}\hline z&1&2 \\ \hline
&-&- \\ \hline
\end{array}$ \ \ \ \ 
$\begin{array}{|c|c|c|c|c|}\hline z&1&2&3&4 \\ \hline
&(6,-6)'&(28,38)*&-&(32,62)* \\ \hline
\end{array}$
\vspace{2mm}

\noindent $p=7$:
\vspace{2mm}

\begin{tabular}{|c|c|c|c|c|c|c|}\hline $z$&1&2&3&4&5&6 \\ \hline
&(2,-46)&(-8,2)&(32,-94)*&(80,290)*&(10,50)'&- \\ \hline
\end{tabular}
\vspace{2mm}

\noindent $p=11$:
\vspace{2mm}

\begin{tabular}{|c|c|c|c|c|c|c|c|c|}\hline $z$&1&2&3&4&5&6&7&8 \\ \hline
&(56,290)'&-&(-16,2)'&(6,26)&(16,98)&(12,114)'&(26,106)&- \\ \hline
\end{tabular}
\vspace{1mm}

\begin{tabular}{|c|c|c|}\hline $z$&9&10 \\ \hline
&(-8,2)&(-36,210)' \\ \hline
\end{tabular}
\vspace{2mm}

\noindent $p=13$:
\vspace{2mm}

\begin{tabular}{|c|c|c|c|c|c|c|c|c|}\hline $z$&1&2&3&4&5&6&7&8 \\ \hline
&(-8,270)'&(20,-106)&(-4,86)&(-204,646)*&(22,-30)&(-160,30)*&(-34,50)&(-16,302) \\ \hline
\end{tabular}
\vspace{1mm}

\begin{tabular}{|c|c|c|c|c|}\hline $z$&9&10&11&12 \\ \hline
&(58,146)&(18,34)&(84,406)&(56,206)' \\ \hline
\end{tabular}
\vspace{2mm}

\noindent $p=17$:
\vspace{2mm}

\begin{tabular}{|c|c|c|c|c|c|c|c|c|}\hline $z$&1&2&3&4&5&6&7&8 \\ \hline
&(256,-322)*&(256,-322)*&(-24,542)&(44,166)&(210,1218)&(24,-178)'&(-100,278)&(22,50) \\ \hline
\end{tabular}
\vspace{1mm}

\begin{tabular}{|c|c|c|c|c|c|c|c|c|}\hline $z$&9&10&11&12&13&14&15&16 \\ \hline
&(-4,70)&(52,470)&(-84,342)'&-&(22,-334)'&(18,258)&(184,974)&(-56,302)' \\ \hline
\end{tabular}
\vspace{2mm}

\subsection{The Case $B*a$}
\normalsize
This is operator nr. 15 from the list \cite{AESZ}:
\[{\theta}^{4}-3\,x \left( 3\,\theta+1 \right)  \left( 3\,\theta+2
 \right)  \left( 7\,{\theta}^{2}+7\,\theta+2 \right) -72\,{x}^{2}
 \left( 3\,\theta+1 \right)  \left( 3\,\theta+2 \right)  \left( 3\,
\theta+4 \right)  \left( 3\,\theta+5 \right)
\]
\scriptsize

$p=3$: \ \ \ \ \ \ \ \ \ \ \ \ \ \ \ \ \ \ \ \ \ $p=5$:
\vspace{2mm}

$\begin{array}{|c|c|c|}\hline z&1&2 \\ \hline
&(2,4)&(8,13) \\ \hline
\end{array}$\ \ \ \
$\begin{array}{|c|c|c|c|c|}\hline z&1&2&3&4 \\ \hline
&(-18,-22)&-&(3,-22)&(6,41) \\ \hline
\end{array}$
\vspace{2mm}

\noindent $p=7$:
\vspace{2mm}

\begin{tabular}{|c|c|c|c|c|c|c|}\hline $z$&1&2&3&4&5&6 \\ \hline
&(-31,-102)&(-13,60)&-&(20,12)&-&- \\ \hline
\end{tabular}
\vspace{2mm}

\noindent $p=11$:
\vspace{2mm}

\begin{tabular}{|c|c|c|c|c|c|c|c|c|}\hline $z$&1&2&3&4&5&6&7&8 \\ \hline
&(36,170)&(-147,422)&(-15,152)&(21,170)&(45,224)&(-24,71)&(-3,-28)&(-72,-478)\\ \hline
\end{tabular}
\vspace{1mm}

\begin{tabular}{|c|c|c|}\hline  $z$ & 9&10\\ \hline
&(51,170)&(-12,8)\\\hline
\end{tabular}

\vspace{2mm}

\noindent $p=13$:
\vspace{2mm}

\begin{tabular}{|c|c|c|c|c|c|c|c|c|}\hline $z$&1&2&3&4&5&6&7&8 \\ \hline
&(23,60)&(20,192)&(-13,72)&(23,330)&(-103,-768)&-&(50,285)&(14,-138) \\ \hline
\end{tabular}
\vspace{1mm}

\begin{tabular}{|c|c|c|c|c|}\hline $z$ & 9 & 10 & 11 & 12 \\ \hline
&(17,144)&(56,228)&-&(-202,618)\\ \hline
\end{tabular}

\vspace{2mm}

\noindent $p=17$:
\vspace{2mm}

\begin{tabular}{|c|c|c|c|c|c|c|c|c|}\hline $z$&1&2&3&4&5&6&7&8 \\ \hline
&(-12,128)&(105,488)&(93,254)&(21,-250)&(-234,-718)&(-60,-25)&(-39,38)&- \\ \hline
\end{tabular}
\vspace{1mm}

\begin{tabular}{|c|c|c|c|c|c|c|c|c|}\hline $z$&9&10&11&12&13&14&15&16 \\ \hline
&(-132,668)&(-414,2522)&(108,362)&(117,524)&(-39,-142)&(-21,488)&-&(15,-196) \\ \hline
\end{tabular}
\vspace{2mm}

\subsection{The Case $C*a$}
\normalsize
This is operator nr. 68 from the list \cite{AESZ}:
\[{\theta}^{4}-4\,x \left( 4\,\theta+1 \right)  \left( 4\,\theta+3
 \right)  \left( 7\,{\theta}^{2}+7\,\theta+2 \right) -128\,{x}^{2}
 \left( 4\,\theta+1 \right)  \left( 4\,\theta+3 \right)  \left( 4\,
\theta+5 \right)  \left( 4\,\theta+7 \right)
 \]
\scriptsize
$p=3$: \ \ \ \ \ \ \ \ \ \ \ \ \ \ \ \ \ \ \ \ $p=5$:
\vspace{2mm}

$\begin{array}{|c|c|c|}\hline z&1&2 \\ \hline
&(2,-2)&- \\ \hline
\end{array}$ \ \ \ \
$\begin{array}{|c|c|c|c|c|}\hline z&1&2&3&4 \\ \hline
&-&(6,6)&(-18,-22)*&(8,38) \\ \hline
\end{array}$
\vspace{2mm}

\noindent $p=7$:
\vspace{2mm}

\begin{tabular}{|c|c|c|c|c|c|c|}\hline $z$&1&2&3&4&5&6 \\ \hline
&(24,-158)*&-&(4,2)'&(22,50)&(-2,66)'&(72,226)* \\ \hline
\end{tabular}
\vspace{2mm}

\noindent $p=11$:
\vspace{2mm}

\begin{tabular}{|c|c|c|c|c|c|c|c|c|}\hline $z$&1&2&3&4&5&6&7&8 \\ \hline
&(10,158)&(124,146)*&-&(-10,-38)&-&(92,-238)*&(28,34)&(-32,122) \\ \hline
\end{tabular}
\vspace{1mm}

\begin{tabular}{|c|c|c|}\hline $z$&9&10 \\ \hline
&-&(14,50) \\ \hline
\end{tabular}
\vspace{2mm}

\noindent $p=13$:
\vspace{2mm}

\begin{tabular}{|c|c|c|c|c|c|c|c|c|}\hline $z$&1&2&3&4&5&6&7&8 \\ \hline
&(232,1038)*&(-4,150)&(-32,62)'&(-46,146)'&(46,126)&(-2,210)&(58,290)&(162,58)* \\ \hline
\end{tabular}
\vspace{1mm}

\begin{tabular}{|c|c|c|c|c|}\hline $z$&9&10&11&12 \\ \hline
&(12,-50)&(64,206)&(-6,86)&(24,262)\\ \hline
\end{tabular}
\vspace{2mm}

\noindent $p=17$:
\vspace{2mm}

\begin{tabular}{|c|c|c|c|c|c|c|c|c|}\hline $z$&1&2&3&4&5&6&7&8 \\ \hline
&(-60,246)&(-30,162)&(52,226)&-&(8,134)&-&(178,962)&- \\ \hline
\end{tabular}
\vspace{1mm}

\begin{tabular}{|c|c|c|c|c|c|c|c|c|}\hline $z$&9&10&11&12&13&14&15&16 \\ \hline
&(404,2342)*&-&-&(-32,-190)&(336,1118)*&(24,-142)&(-24,254)&(66,506) \\ \hline
\end{tabular}
\vspace{2mm}

\subsection{The Case $D*a$ }
\normalsize
This is operator nr. 62 from the list \cite{AESZ}:
\[{\theta}^{4}-12\,x \left( 6\,\theta+1 \right)  \left( 6\,\theta+5
 \right)  \left( 7\,{\theta}^{2}+7\,\theta+2 \right) -1152\,{x}^{2}
 \left( 6\,\theta+1 \right)  \left( 6\,\theta+5 \right)  \left( 6\,
\theta+7 \right)  \left( 6\,\theta+11 \right)
\]
\scriptsize
$p=3$: \ \ \ \ \ \ \ \ \ \ \ \ \ \ \ \ \ \ \ \ \ \ $p=5$:
\vspace{2mm}

$\begin{array}{|c|c|c|}\hline z&1&2 \\ \hline
&(2,4)&(8,13) \\ \hline
\end{array}$\ \ \ \
$\begin{array}{|c|c|c|c|c|}\hline z&1&2&3&4 \\ \hline
&(34,74)*&(29,44)*&-&- \\ \hline
\end{array}$
\vspace{2mm}

\noindent $p=7$:
\vspace{2mm}

\begin{tabular}{|c|c|c|c|c|c|c|}\hline $z$&1&2&3&4&5&6 \\ \hline
&(5,-4)&(4,-40)&(59,122)*&(65,170)*&(22,92)&- \\ \hline
\end{tabular}
\vspace{2mm}

\noindent $p=11$:
\vspace{2mm}

\begin{tabular}{|c|c|c|c|c|c|c|c|c|}\hline $z$&1&2&3&4&5&6&7&8 \\ \hline
&-&(12,96)&(-9,-46)&(25,14)&(59,296)&(-160,578)*&(-115,38)*&(29,184) \\ \hline
\end{tabular}
\vspace{1mm}

\begin{tabular}{|c|c|c|}\hline $z$&9&10 \\ \hline
&(8,-142)&(-24,15) \\ \hline
\end{tabular}
\vspace{2mm}

\noindent $p=13$:
\vspace{2mm}

\begin{tabular}{|c|c|c|c|c|c|c|c|c|}\hline $z$&1&2&3&4&5&6&7&8 \\ \hline
&(67,276)&(56,374)&(5,-100)'&(-138,-278)*&(4,-12)&(-193,492)*&(38,350)&(47,188) \\ \hline
\end{tabular}
\vspace{1mm}

\begin{tabular}{|c|c|c|c|c|}\hline $z$&9&10&11&12 \\ \hline
&(-23,8)'&(-3,322)'&(-36,199)&(8,-160) \\ \hline
\end{tabular}
\vspace{2mm}

\noindent $p=17$:
\vspace{2mm}

\begin{tabular}{|c|c|c|c|c|c|c|c|c|}\hline $z$&1&2&3&4&5&6&7&8 \\ \hline
&(67,284)&(18,-79)&(-131,728)&(45,218)&(19,-388)'&(80,490)&(262,-214)*&(72,164) \\ \hline
\end{tabular}
\vspace{1mm}

\begin{tabular}{|c|c|c|c|c|c|c|c|c|}\hline $z$&9&10&11&12&13&14&15&16 \\ \hline
&(-150,822)&(55,-70)'&(160,863)&(250,-430)*&(-15,150)&(11,-278)&(141,768)&(-16,56) \\ \hline
\end{tabular}
\vspace{1mm}

\subsection{The Case $A*b$}
\normalsize
This is operator nr. 25 from the list \cite{AESZ}:
\[{\theta}^{4}-4\,x \left( 2\,\theta+1 \right) ^{2} \left( 11\,{\theta}^
{2}+11\,\theta+3 \right) -16\,{x}^{2} \left( 2\,\theta+1 \right) ^{2}
 \left( 2\,\theta+3 \right) ^{2}
\]
\scriptsize
$p=3$: \ \ \ \ \ \ \ \ \ \ \ \ \ \ \ \ \ \ \ \ \ $p=5$:
\vspace{2mm}

$\begin{array}{|c|c|c|}\hline z&1&2 \\ \hline
&-&(5,14) \\ \hline
\end{array}$ \ \ \ \
$\begin{array}{|c|c|c|c|c|}\hline z&1&2&3&4 \\ \hline
&(13,16)&-&(2,26)'&(-3,16) \\ \hline
\end{array}$
\vspace{2mm}

\noindent $p=7$:
\vspace{2mm}

\begin{tabular}{|c|c|c|c|c|c|c|}\hline $z$&1&2&3&4&5&6 \\ \hline
&(10,50)'&(25,74)&-&(-10,82)&(10,-30)&(-15,26) \\ \hline
\end{tabular}
\vspace{2mm}
\newpage
\noindent $p=11$:
\vspace{2mm}

\begin{tabular}{|c|c|c|c|c|c|c|c|c|}\hline $z$&1&2&3&4&5&6&7&8 \\ \hline
&(39,262)&(112,2)*&(2,58)'&(-26,42)&(15,166)&(10,-134)&(39,142)&(47,78) \\ \hline
\end{tabular}
\vspace{1mm}

\begin{tabular}{|c|c|c|}\hline $z$&9&10 \\ \hline
&(152,482)*&(-26,122) \\ \hline
\end{tabular}
\vspace{2mm}

\noindent $p=13$:
\vspace{2mm}

\begin{tabular}{|c|c|c|c|c|c|c|c|c|}\hline $z$&1&2&3&4&5&6&7&8 \\ \hline
&(-60,166)'&-&-&-&-&(-30,90)&(20,214)&(35,120) \\ \hline
\end{tabular}
\vspace{1mm}

\begin{tabular}{|c|c|c|c|c|}\hline $z$&9&10&11&12 \\ \hline
&(50,98)&(60,246)&(35,-40)&(15,12) \\ \hline
\end{tabular}
\vspace{2mm}

\noindent $p=17$:
\vspace{2mm}

\begin{tabular}{|c|c|c|c|c|c|c|c|c|}\hline $z$&1&2&3&4&5&6&7&8 \\ \hline
&-&(115,744)&(20,86)'&(-10,50)&(-15,368)&(-25,60)&(140,790)'&(35,-56) \\ \hline
\end{tabular}
\vspace{1mm}

\begin{tabular}{|c|c|c|c|c|c|c|c|c|}\hline $z$&9&10&11&12&13&14&15&16 \\ \hline
&-&(-15,208)&(20,-394)&(-20,-330)&-&(55,540)&(75,632)&(60,134) \\ \hline
\end{tabular}
\vspace{1mm}

\subsection{The Case $B*b$}
\normalsize
This is operator nr. 24 from the list \cite{AESZ}:
\[{\theta}^{4}-3\,x \left( 3\,\theta+1 \right)  \left( 3\,\theta+2
 \right)  \left( 11\,{\theta}^{2}+11\,\theta+3 \right) -9\,{x}^{2}
 \left( 3\,\theta+1 \right)  \left( 3\,\theta+2 \right)  \left( 3\,
\theta+4 \right)  \left( 3\,\theta+5 \right)
\]
\scriptsize
$p=3$:\ \ \ \ \ \ \ \ \ \ \ \ \ \ \ \ \ \ \ \ \ \ \ $p=5$:
\vspace{2mm}

$\begin{array}{|c|c|c|}\hline z&1&2 \\ \hline
&(5,7)&(5,19) \\ \hline
\end{array}$\ \ \ \
$\begin{array}{|c|c|c|c|c|}\hline z&1&2&3&4 \\ \hline
&-&(8,-4)&-&(7,-4) \\ \hline
\end{array}$
\vspace{2mm}

\noindent $p=7$:
\vspace{2mm}

\begin{tabular}{|c|c|c|c|c|c|c|}\hline $z$&1&2&3&4&5&6 \\ \hline
&-&(25,113)&(25,86)&-&(-15,-11)&(15,25) \\ \hline
\end{tabular}
\vspace{2mm}

\noindent $p=11$:
\vspace{2mm}

\begin{tabular}{|c|c|c|c|c|c|c|c|c|}\hline $z$&1&2&3&4&5&6&7&8 \\ \hline
&-&(105,-82)&(-29,152)&(10,-127)&(-3,62)&(37,197)&(15,188)&(36,107) \\ \hline
\end{tabular}
\vspace{1mm}

\begin{tabular}{|c|c|c|}\hline $z$&9&10 \\ \hline
&(150,458)&- \\ \hline
\end{tabular}
\vspace{2mm}

\noindent $p=13$:
\vspace{2mm}

\begin{tabular}{|c|c|c|c|c|c|c|c|c|}\hline $z$&1&2&3&4&5&6&7&8 \\ \hline
&(90,319)&(15,112)&-&(-35,-4)&-&(45,142)&-&(-5,-151) \\ \hline
\end{tabular}
\vspace{1mm}

\begin{tabular}{|c|c|c|c|c|}\hline $z$&9&10&11&12 \\ \hline
&(20,210)&(35,252)&(-85,447)&(-45,49) \\ \hline
\end{tabular}
\vspace{2mm}

\noindent $p=17$:
\vspace{2mm}

\begin{tabular}{|c|c|c|c|c|c|c|c|c|}\hline $z$&1&2&3&4&5&6&7&8 \\ \hline
&(50,-115)&(-30,524)&(10,362)&(-10,83)&(65,470)&(165,947)&(30,362)&(10,362) \\ \hline
\end{tabular}
\vspace{1mm}

\begin{tabular}{|c|c|c|c|c|c|c|c|c|}\hline $z$&9&10&11&12&13&14&15&16 \\ \hline
&(80,407)&(45,83)&-&(110,461)&(-120,569)&(-25,-196)&(-40,38)&(-50,56) \\ \hline
\end{tabular}
\vspace{1mm}

\subsection{The Case $C*b$}
\normalsize
This is operator nr. 51 from the list \cite{AESZ}:
\[{\theta}^{4}-4\,x \left( 4\,\theta+1 \right)  \left( 4\,\theta+3
 \right)  \left( 11\,{\theta}^{2}+11\,\theta+3 \right) -16\,{x}^{2}
 \left( 4\,\theta+1 \right)  \left( 4\,\theta+3 \right)  \left( 4\,
\theta+5 \right)  \left( 4\,\theta+7 \right)
\]
\scriptsize

$p=3$: \ \ \ \ \ \ \ \ \ \ \ \ \ \ \ \ \ \ \ \ \ $p=5$:
\vspace{2mm}

$\begin{array}{|c|c|c|}\hline z&1&2 \\ \hline
&(5,14)&(5,2) \\ \hline
\end{array}$\ \ \ \
$\begin{array}{|c|c|c|c|c|}\hline z&1&2&3&4 \\ \hline
&(3,-4)'&(12,46)&-&- \\ \hline
\end{array}$
\vspace{2mm}

\noindent $p=7$:
\vspace{2mm}

\begin{tabular}{|c|c|c|c|c|c|c|}\hline $z$&1&2&3&4&5&6 \\ \hline
&(40,122)&(-10,18)&-&(-5,90)&(15,26)&- \\ \hline
\end{tabular}
\vspace{2mm}

\noindent $p=11$:
\vspace{2mm}

\begin{tabular}{|c|c|c|c|c|c|c|c|c|}\hline $z$&1&2&3&4&5&6&7&8 \\ \hline
&(24,130)&(39,162)'&-&(-5,-74)&(-64,-574)*&(-144,386)*&(30,206)&(-4,162) \\ \hline
\end{tabular}
\vspace{1mm}

\begin{tabular}{|c|c|c|}\hline $z$&9&10 \\ \hline
&(-26,122)&(19,130) \\ \hline
\end{tabular}
\vspace{2mm}

\noindent $p=13$:
\vspace{2mm}

\begin{tabular}{|c|c|c|c|c|c|c|c|c|}\hline $z$&1&2&3&4&5&6&7&8 \\ \hline
&(80,430)&(75,282)&(-15,96)&(45,228)&(-30,-38)&(-5,-190)&(30,166)&(-80,282) \\ \hline
\end{tabular}
\vspace{1mm}

\begin{tabular}{|c|c|c|c|c|}\hline $z$&9&10&11&12 \\ \hline
&(30,202)&(10,202)&(30,122)&- \\ \hline
\end{tabular}
\vspace{2mm}

\noindent $p=17$:
\vspace{2mm}

\begin{tabular}{|c|c|c|c|c|c|c|c|c|}\hline $z$&1&2&3&4&5&6&7&8 \\ \hline
&(10,522)&(35,292)&(70,626)&(90,554)&(-70,382)&(50,-110)'&(90,326)&(-25,188) \\ \hline
\end{tabular}
\vspace{1mm}

\begin{tabular}{|c|c|c|c|c|c|c|c|c|}\hline $z$&9&10&11&12&13&14&15&16 \\ \hline
&-&(65,514)&(15,-150)'&(65,450)&(-50,162)&(115,410)'&(15,124)&(-100,326) \\ \hline
\end{tabular}
\vspace{1mm}

\subsection{The Case $D*b$}
\normalsize
This is operator nr. 63 from the list \cite{AESZ}:
\[{\theta}^{4}-12\,x \left( 6\,\theta+1 \right)  \left( 6\,\theta+5
 \right)  \left( 11\,{\theta}^{2}+11\,\theta+3 \right) -144\,{x}^{2}
 \left( 6\,\theta+1 \right)  \left( 6\,\theta+5 \right)  \left( 6\,
\theta+7 \right)  \left( 6\,\theta+11 \right)
\]
\scriptsize

$p=3$:\ \ \ \ \ \ \ \ \ \ \ \ \ \ \ \ \ \ \ \ \ \ $p=5$:
\vspace{2mm}

$\begin{array}{|c|c|c|}\hline z&1&2 \\ \hline
&(5,7)&(5,19) \\ \hline
\end{array}$\ \ \ \
$\begin{array}{|c|c|c|c|c|}\hline z&1&2&3&4 \\ \hline
&-&(24,76)&(4,1)&(-1,-4) \\ \hline
\end{array}$
\vspace{2mm}

\noindent $p=7$:
\vspace{2mm}

\begin{tabular}{|c|c|c|c|c|c|c|}\hline $z$&1&2&3&4&5&6 \\ \hline
&(15,47)&-&(5,31)&-&(-5,62)&(25,95) \\ \hline
\end{tabular}
\vspace{2mm}

\noindent $p=11$:
\vspace{2mm}

\begin{tabular}{|c|c|c|c|c|c|c|c|c|}\hline $z$&1&2&3&4&5&6&7&8 \\ \hline
&(39,142)&(13,137)&(-2,87)&(104,-94)*&(23,4)&(8,129)&(169,686)*&- \\ \hline
\end{tabular}
\vspace{1mm}

\begin{tabular}{|c|c|c|}\hline $z$&9&10 \\ \hline
&-&(-41,157) \\ \hline
\end{tabular}
\vspace{2mm}

\noindent $p=13$:
\vspace{2mm}

\begin{tabular}{|c|c|c|c|c|c|c|c|c|}\hline $z$&1&2&3&4&5&6&7&8 \\ \hline
&(15,139)&(85,410)'&(40,86)&(75,275)&(-5,-268)&-&(-55,355)&(-25,189) \\ \hline
\end{tabular}
\vspace{1mm}

\begin{tabular}{|c|c|c|c|c|}\hline $z$&9&10&11&12 \\ \hline
&(20,293)&(15,-120)&(-40,305)&(-15,-180) \\ \hline
\end{tabular}
\vspace{2mm}

\noindent $p=17$:
\vspace{2mm}

\begin{tabular}{|c|c|c|c|c|c|c|c|c|}\hline $z$&1&2&3&4&5&6&7&8 \\ \hline
&(30,88)&(-10,206)&(15,236)&(20,111)&(90,239)&(140,749)&-&(10,-231) \\ \hline
\end{tabular}
\vspace{1mm}

\begin{tabular}{|c|c|c|c|c|c|c|c|c|}\hline $z$&9&10&11&12&13&14&15&16 \\ \hline
&-&(-30,410)&(5,-41)&(105,698)&-&(-10,542)&(-140,684)&(50,-137) \\ \hline
\end{tabular}
\vspace{1mm}

\subsection{The Case $A*c$ }
\normalsize
This is operator nr.58 from the list \cite{AESZ}:
\[
{\theta}^{4}-4\,x \left( 2\,\theta+1 \right) ^{2} \left( 10\,{\theta}^
{2}+10\,\theta+3 \right) +144\,{x}^{2} \left( 2\,\theta+1 \right) ^{2}
 \left( 2\,\theta+3 \right) ^{2}
\]
\scriptsize
$p=3$:\ \ \ \ \ \ \ \ \ \ \ \ \ \ \ \ \  \ \ \ \ $p=5$:
\vspace{2mm}

$\begin{array}{|c|c|c|}\hline z&1&2 \\ \hline
&(8,2)*&- \\ \hline
\end{array}$\ \ \ \
$\begin{array}{|c|c|c|c|c|}\hline z&1&2&3&4 \\ \hline
&(-28,38)*&-&(-2,-14)&(16,-34)* \\ \hline
\end{array}$
\vspace{2mm}

\noindent $p=7$:
\vspace{2mm}

\begin{tabular}{|c|c|c|c|c|c|c|}\hline $z$&1&2&3&4&5&6 \\ \hline
&(12,2)&(32,-94)*&(26,50)&(80,290)*&-&(12,82) \\ \hline
\end{tabular}
\vspace{2mm}

\noindent $p=11$:
\vspace{2mm}

\begin{tabular}{|c|c|c|c|c|c|c|c|c|}\hline $z$&1&2&3&4&5&6&7&8 \\ \hline
&(-160,578)*&(-60,290)&(-12,-78)&(4,-14)&(20,178)&(14,122)'&(-46,170)&(-4,-14) \\ \hline
\end{tabular}
\vspace{1mm}

\begin{tabular}{|c|c|c|}\hline $z$&9&10 \\ \hline
&-&(36,98) \\ \hline
\end{tabular}
\vspace{2mm}

\noindent $p=13$:
\vspace{2mm}

\begin{tabular}{|c|c|c|c|c|c|c|c|c|}\hline $z$&1&2&3&4&5&6&7&8 \\ \hline
&(-108,-698)*&(-14,-70)&(16,126)&(32,158)&(8,62)&(42,202)'&(16,62)&(42,10) \\ \hline
\end{tabular}
\vspace{1mm}

\begin{tabular}{|c|c|c|c|c|}\hline $z$&9&10&11&12 \\ \hline
&(-204,646)*&(16,126)&(2,314)&(-16,254)' \\ \hline
\end{tabular}
\vspace{2mm}

\noindent $p=17$:
\vspace{2mm}

\begin{tabular}{|c|c|c|c|c|c|c|c|c|}\hline $z$&1&2&3&4&5&6&7&8 \\ \hline
&(-76,278)&(-8,-178)&(-134,562)'&(8,302)&(-24,-178)'&(-142,706)&(-32,110)&(168,942) \\ \hline
\end{tabular}
\vspace{1mm}

\begin{tabular}{|c|c|c|c|c|c|c|c|c|}\hline $z$&9&10&11&12&13&14&15&16 \\ \hline
&-&(76,278)&(66,2)&(-38,178)&(-12,-234)&-&(224,-898)*&(-356,1478)* \\ \hline
\end{tabular}
\vspace{1mm}

\subsection{The Case $B*c$}
\normalsize
This is operator nr.70 from the list \cite{AESZ}:
\[
{\theta}^{4}-3\,x \left( 3\,\theta+1 \right)  \left( 3\,\theta+2
 \right)  \left( 10\,{\theta}^{2}+10\,\theta+3 \right) +81\,{x}^{2}
 \left( 3\,\theta+1 \right)  \left( 3\,\theta+2 \right)  \left( 3\,
\theta+4 \right)  \left( 3\,\theta+5 \right)
\]
\scriptsize
$p=3$: \ \ \ \ \ \ \ \ \ \ \ \ \ \ \ \ \ \ \ \ \ \ \ \ \ \ \ $p=5$:
\vspace{2mm}

$\begin{array}{|c|c|c|}\hline z&1&2 \\ \hline
&(5,10)'&(-4,-2)' \\ \hline
\end{array}$\ \ \ \
$\begin{array}{|c|c|c|c|c|}\hline z&1&2&3&4 \\ \hline
&(-9,-4)&(-27,32)*&(11,-16)'&(-3,32) \\ \hline
\end{array}$
\vspace{2mm}

\noindent $p=7$:
\vspace{2mm}

\begin{tabular}{|c|c|c|c|c|c|c|}\hline $z$&1&2&3&4&5&6 \\ \hline
&(17,54)'&(2,30)&(-46,18)*&-&(11,-24)&(-31,-102)* \\ \hline
\end{tabular}
\vspace{2mm}

\noindent $p=11$:
\vspace{2mm}

\begin{tabular}{|c|c|c|c|c|c|c|c|c|}\hline $z$&1&2&3&4&5&6&7&8 \\ \hline
&(-144,386)*&(18,89)'&(3,-100)&-&-&(-6,26)&-&(-72,350) \\ \hline
\end{tabular}
\vspace{1mm}

\begin{tabular}{|c|c|c|}\hline $z$&9&10 \\ \hline
&(147,422)*&(27,62)' \\ \hline
\end{tabular}
\vspace{2mm}

\newpage
\noindent $p=13$:
\vspace{2mm}

\begin{tabular}{|c|c|c|c|c|c|c|c|c|}\hline $z$&1&2&3&4&5&6&7&8 \\ \hline
&(-202,618)*&(62,198)&(-190,450)*&(-34,147)&(5,150)'&-&(20,30)&- \\ \hline
\end{tabular}
\vspace{1mm}

\begin{tabular}{|c|c|c|c|c|}\hline $z$&9&10&11&12 \\ \hline
&(20,30)&(-31,203)&-&(41,240) \\ \hline
\end{tabular}
\vspace{2mm}

\noindent $p=17$:
\vspace{2mm}

\begin{tabular}{|c|c|c|c|c|c|c|c|c|}\hline $z$&1&2&3&4&5&6&7&8 \\ \hline
&-&(174,947)'&(-33,20)&(39,326)&(-16,57)&(6,362)&(-180,-1690)*&(-18,2) \\ \hline
\end{tabular}
\vspace{1mm}

\begin{tabular}{|c|c|c|c|c|c|c|c|c|}\hline $z$&9&10&11&12&13&14&15&16 \\ \hline
&(-18,-358)&(-63,200)&-&(234,-718)*&(-39,-214)&(-81,92)&(-135,776)'&(-144,866) \\ \hline
\end{tabular}
\vspace{1mm}

\subsection{The Case $C*c$}
\normalsize
This is operator nr. 69 from the list \cite{AESZ}:
\[
{\theta}^{4}-4\,x \left( 4\,\theta+1 \right)  \left( 4\,\theta+3
 \right)  \left( 10\,{\theta}^{2}+10\,\theta+3 \right) +144\,{x}^{2}
 \left( 4\,\theta+1 \right)  \left( 4\,\theta+3 \right)  \left( 4\,
\theta+5 \right)  \left( 4\,\theta+7 \right)
\]
\scriptsize
$p=3$:\ \ \ \ \ \ \ \ \ \ \ \ \ \ \ \ \ \ \ \ \ \ \ \ \ \  \ \ \ \ \  \ \ $p=5$:
\vspace{2mm}

$\begin{array}{|c|c|c|}\hline z&1&2 \\ \hline
&(-4,-14)*&(-4,10) \\ \hline
\end{array}$\ \ \ \
$\begin{array}{|c|c|c|c|c|}\hline z&1&2&3&4 \\ \hline
&(-32,62)*&(-8,2)'&(-4,26)&- \\ \hline
\end{array}$
\vspace{2mm}

\noindent $p=7$:
\vspace{2mm}

\begin{tabular}{|c|c|c|c|c|c|c|}\hline $z$&1&2&3&4&5&6 \\ \hline
&(40,-30)*&-&(-4,-54)&(44,2)*&(36,118)&- \\ \hline
\end{tabular}
\vspace{2mm}

\noindent $p=11$:
\vspace{2mm}

\begin{tabular}{|c|c|c|c|c|c|c|c|c|}\hline $z$&1&2&3&4&5&6&7&8 \\ \hline
&(32,130)&(16,2)'&(72,-478)*&(-20,-46)&(-172,722)*&(12,54)&(-40,218)&(-20,182) \\ \hline
\end{tabular}
\vspace{1mm}

\begin{tabular}{|c|c|c|}\hline $z$&9&10 \\ \hline
&(28,82)&(-28,50) \\ \hline
\end{tabular}
\vspace{2mm}

\noindent $p=13$:
\vspace{2mm}

\begin{tabular}{|c|c|c|c|c|c|c|c|c|}\hline $z$&1&2&3&4&5&6&7&8 \\ \hline
&(-20,-138)&(4,218)&-&(40,206)'&(40,2)&(72,290)&-&- \\ \hline
\end{tabular}
\vspace{1mm}

\begin{tabular}{|c|c|c|c|c|}\hline $z$&9&10&11&12 \\ \hline
&(-12,70)&(140,-250)*&(60,334)&(132,-362)* \\ \hline
\end{tabular}
\vspace{2mm}

\noindent $p=17$:
\vspace{2mm}

\begin{tabular}{|c|c|c|c|c|c|c|c|c|}\hline $z$&1&2&3&4&5&6&7&8 \\ \hline
&(-24,-82)&(-72,110)&(-12,26)&(-276,38)*&(-76,122)&(148,734)&(88,218)&(316,758)* \\ \hline
\end{tabular}
\vspace{1mm}

\begin{tabular}{|c|c|c|c|c|c|c|c|c|}\hline $z$&9&10&11&12&13&14&15&16 \\ \hline
&(-28,-58)&(-176,962)&(-112,386)'&(-28,470)&(-120,462)&(24,210)'&(-4,-266)&(64,382)' \\ \hline
\end{tabular}
\vspace{1mm}

\subsection{The Case $D*c$}
\normalsize
This is operator nr. 64 from the list \cite{AESZ}:
\[
{\theta}^{4}-12\,x \left( 6\,\theta+1 \right)  \left( 6\,\theta+5
 \right)  \left( 10\,{\theta}^{2}+10\,\theta+3 \right) +1296\,{x}^{2}
 \left( 6\,\theta+1 \right)  \left( 6\,\theta+5 \right)  \left( 6\,
\theta+7 \right)  \left( 6\,\theta+11 \right)
\]
\scriptsize
$p=3$:\ \ \ \ \ \ \ \ \ \ \ \ \ \ \ \ \ \ \ \ \ \ \ \ \ \ \ \ \ \ $p=5$:
\vspace{5mm}

$\begin{array}{|c|c|c|}\hline z&1&2 \\ \hline
&(5,10)'&(-4,-2)' \\ \hline
\end{array}$\ \ \ \
$\begin{array}{|c|c|c|c|c|}\hline z&1&2&3&4 \\ \hline
&-&(19,-16)*&(-31,56)*&- \\ \hline
\end{array}$
\vspace{2mm}

\noindent $p=7$:
\vspace{2mm}

\begin{tabular}{|c|c|c|c|c|c|c|}\hline $z$&1&2&3&4&5&6 \\ \hline
&(-6,-50)&(31,128)'&(47,26)*&-&(86,338)*&- \\ \hline
\end{tabular}
\vspace{2mm}

\noindent $p=11$:
\vspace{2mm}

\begin{tabular}{|c|c|c|c|c|c|c|c|c|}\hline $z$&1&2&3&4&5&6&7&8 \\ \hline
&(-49,238)&(-75,350)&(31,76)&(115,38)*&(-21,60)&(8,-98)&(-18,-7)&- \\ \hline
\end{tabular}
\vspace{1mm}

\begin{tabular}{|c|c|c|}\hline $z$&9&10 \\ \hline
&(-136,290)*&(14,122)' \\ \hline
\end{tabular}
\vspace{2mm}

\noindent $p=13$:
\vspace{2mm}

\begin{tabular}{|c|c|c|c|c|c|c|c|c|}\hline $z$&1&2&3&4&5&6&7&8 \\ \hline
&(-198,562)*&-&(-44,222)&(-31,8)&-&(75,310)&(25,140)&(45,160) \\ \hline
\end{tabular}
\vspace{1mm}

\begin{tabular}{|c|c|c|c|c|}\hline $z$&9&10&11&12 \\ \hline
&(-138,-278)*&(22,75)&(44,254)&(-7,-4) \\ \hline
\end{tabular}
\vspace{2mm}

\noindent $p=17$:
\vspace{2mm}

\begin{tabular}{|c|c|c|c|c|c|c|c|c|}\hline $z$&1&2&3&4&5&6&7&8 \\ \hline
&(16,-94)&(121,520)&(-111,444)&-&(-362,1586)*&-&(79,488)&(-2,250) \\ \hline
\end{tabular}
\vspace{1mm}

\begin{tabular}{|c|c|c|c|c|c|c|c|c|}\hline $z$&9&10&11&12&13&14&15&16 \\ \hline
&-&(236,-682)*&(-6,-342)&(95,392)&(63,254)&-&(-162,851)&(-83,368) \\ \hline
\end{tabular}
\vspace{1mm}

\subsection{The Case $A*d$ }
\normalsize
This is operator nr. 36 from the list \cite{AESZ}:
\[
{\theta}^{4}-16\,x \left( 2\,\theta+1 \right) ^{2} \left( 3\,{\theta}^
{2}+3\,\theta+1 \right) +512\,{x}^{2} \left( 2\,\theta+1 \right) ^{2}
 \left( 2\,\theta+3 \right) ^{2}
\]
\scriptsize
$p=3$: \ \ \ \ \ \ \ \ \ \ \ \ \ \ \ \ \ \ \ \ \ $p=5$:
\vspace{2mm}

$\begin{array}{|c|c|c|}\hline z&1&2 \\ \hline
&(4,-14)&- \\ \hline
\end{array}$\ \ \ \
$\begin{array}{|c|c|c|c|c|}\hline z&1&2&3&4 \\ \hline
&(8,46)&(-8,-82)&-&- \\ \hline
\end{array}$
\vspace{2mm}

\noindent $p=7$:
\vspace{2mm}

\begin{tabular}{|c|c|c|c|c|c|c|}\hline $z$&1&2&3&4&5&6 \\ \hline
&(40,-30)&(-8,-30)&(-12,34)&-&-&- \\ \hline
\end{tabular}
\vspace{2mm}

\noindent $p=11$:
\vspace{2mm}

\begin{tabular}{|c|c|c|c|c|c|c|c|c|}\hline $z$&1&2&3&4&5&6&7&8 \\ \hline
&(-80,322)&(-8,162)&(-28,146)&(-20,82)&(172,722)&-&(16,-30)&- \\ \hline
\end{tabular}
\vspace{1mm}

\begin{tabular}{|c|c|c|}\hline $z$&9&10 \\ \hline
&(24,-62)&(-4,-142) \\ \hline
\end{tabular}
\vspace{2mm}

\noindent $p=13$:
\vspace{2mm}

\begin{tabular}{|c|c|c|c|c|c|c|c|c|}\hline $z$&1&2&3&4&5&6&7&8 \\ \hline
&(-36,86)&(28,118)&(56,270)&(36,230)&(-48,254)&(-200,590)&(72,398)&- \\ \hline
\end{tabular}
\vspace{1mm}

\begin{tabular}{|c|c|c|c|c|}\hline $z$&9&10&11&12 \\ \hline
&-&(-18,8)&(60,214)&(-132,-362) \\ \hline
\end{tabular}
\vspace{2mm}

\noindent $p=17$:
\vspace{2mm}

\begin{tabular}{|c|c|c|c|c|c|c|c|c|}\hline $z$&1&2&3&4&5&6&7&8 \\ \hline
&(44,-90)&(-212,-1114)&(-76,598)&(-276,38)&-&(28,326)&-&(-84,422) \\ \hline
\end{tabular}
\vspace{1mm}

\begin{tabular}{|c|c|c|c|c|c|c|c|c|}\hline $z$&9&10&11&12&13&14&15&16 \\ \hline
&(-4,6)&(112,606)&(-16,-162)&(-8,-50)&(-44,598)&(44,-42)&(20,470)&(124,774) \\ \hline
\end{tabular}
\vspace{1mm}

\subsection{The Case $B*d$}
\normalsize
This is operator nr. 48 from the list \cite{AESZ}:
\[
{\theta}^{4}-12\,x \left( 3\,\theta+1 \right)  \left( 3\,\theta+2
 \right)  \left( 3\,{\theta}^{2}+3\,\theta+1 \right) +288\,{x}^{2}
 \left( 3\,\theta+1 \right)  \left( 3\,\theta+2 \right)  \left( 3\,
\theta+4 \right)  \left( 3\,\theta+5 \right)
\]
\scriptsize
$p=3$:\ \ \ \ \ \ \ \ \ \ \ \ \ \ \ \ \ \ \ \  \ \ \ \ \ \ \ \ \ \ \ $p=5$:
\vspace{2mm}

$\begin{array}{|c|c|c|}\hline z&1&2 \\ \hline
&(-1,-8)&(-7,16) \\ \hline
\end{array}$\ \ \ \
$\begin{array}{|c|c|c|c|c|}\hline z&1&2&3&4 \\ \hline
&-&(21,-4)*&-&- \\ \hline
\end{array}$
\vspace{2mm}

\noindent $p=7$:
\vspace{2mm}

\begin{tabular}{|c|c|c|c|c|c|c|}\hline $z$&1&2&3&4&5&6 \\ \hline
&(-5,32)&(-11,32)&(-5,38)&(-8,62)&(-55,90)*&(36,-62)* \\ \hline
\end{tabular}
\vspace{2mm}

\noindent $p=11$:
\vspace{2mm}

\begin{tabular}{|c|c|c|c|c|c|c|c|c|}\hline $z$&1&2&3&4&5&6&7&8 \\ \hline
&-&(-29,152)&(37,80)&(-89,386)&(69,-514)*&(8,-145)&(50,98)'&- \\ \hline
\end{tabular}
\vspace{1mm}

\begin{tabular}{|c|c|c|}\hline $z$&9&10 \\ \hline
&(-40,170)&(-1,98) \\ \hline
\end{tabular}
\vspace{2mm}

\noindent $p=13$:
\vspace{2mm}

\begin{tabular}{|c|c|c|c|c|c|c|c|c|}\hline $z$&1&2&3&4&5&6&7&8 \\ \hline
&(36,49)&(21,-44)&(18,322)'&-&(-112,-642)*&(58,98)&-&(-21,334) \\ \hline
\end{tabular}
\vspace{1mm}

\begin{tabular}{|c|c|c|c|c|}\hline $z$&9&10&11&12 \\ \hline
&(27,-56)&(-154,-54)*&(33,166)&(-24,106) \\ \hline
\end{tabular}
\vspace{2mm}

\noindent $p=17$:
\vspace{2mm}

\begin{tabular}{|c|c|c|c|c|c|c|c|c|}\hline $z$&1&2&3&4&5&6&7&8 \\ \hline
&(88,614)&(-32,326)&(234,-718)*&(-11,128)&(-14,-286)&(109,362)'&(-35,146)&(105,308) \\ \hline
\end{tabular}
\vspace{1mm}

\begin{tabular}{|c|c|c|c|c|c|c|c|c|}\hline $z$&9&10&11&12&13&14&15&16 \\ \hline
&(15,20)&-&-&(18,155)&(88,569)&(-5,506)&(-71,452)&(-20,-250) \\ \hline
\end{tabular}
\vspace{1mm}

\subsection{The Case $C*d$}
\normalsize
This is operator nr. 38 from the list \cite{AESZ}:
\[
{\theta}^{4}-16\,x \left( 4\,\theta+1 \right)  \left( 4\,\theta+3
 \right)  \left( 3\,{\theta}^{2}+3\,\theta+1 \right) +512\,{x}^{2}
 \left( 4\,\theta+1 \right)  \left( 4\,\theta+3 \right)  \left( 4\,
\theta+5 \right)  \left( 4\,\theta+7 \right)
\]
\scriptsize
$p=3$: \ \ \ \ \ \ \ \ \ \ \ \ \ \ \ \ \ \ \ \ \ \ \ \ \ \ \ \ \ $p=5$:
\vspace{2mm}

$\begin{array}{|c|c|c|}\hline z&1&2 \\ \hline
&(-10,10)&(2,-22) \\ \hline
\end{array}$\ \ \ \
$\begin{array}{|c|c|c|c|c|}\hline z&1&2&3&4 \\ \hline
&(36,86)&-&-&- \\ \hline
\end{array}$
\vspace{2mm}

\noindent $p=7$:
\vspace{2mm}

\begin{tabular}{|c|c|c|c|c|c|c|}\hline $z$&1&2&3&4&5&6 \\ \hline
&-&(36,-62)&(-12,-2)&(-4,66)&(-10,10)&(-2,26) \\ \hline
\end{tabular}
\vspace{2mm}

\noindent $p=11$:
\vspace{2mm}

\begin{tabular}{|c|c|c|c|c|c|c|c|c|}\hline $z$&1&2&3&4&5&6&7&8 \\ \hline
&(10,122)&(150,458)&(12,-78)&(-118,74)&(-64,306)&(20,146)&(-42,122)&- \\ \hline
\end{tabular}
\vspace{1mm}

\begin{tabular}{|c|c|c|}\hline $z$&9&10 \\ \hline
&(-98,434)&(-4,-30) \\ \hline
\end{tabular}
\vspace{2mm}

\noindent $p=13$:
\vspace{2mm}

\begin{tabular}{|c|c|c|c|c|c|c|c|c|}\hline $z$&1&2&3&4&5&6&7&8 \\ \hline
&(-16,126)&(32,158)&(236,1094)&-&-&(-16,158)&-&- \\ \hline
\end{tabular}
\vspace{1mm}

\begin{tabular}{|c|c|c|c|c|}\hline $z$&9&10&11&12 \\ \hline
&(-14,-86)&(12,54)&(2,42)&(62,346) \\ \hline
\end{tabular}
\vspace{2mm}

\noindent $p=17$:
\vspace{2mm}

\begin{tabular}{|c|c|c|c|c|c|c|c|c|}\hline $z$&1&2&3&4&5&6&7&8 \\ \hline
&(-240,-610)&(96,382)&(62,314)&(-24,14)&(8,78)&(-24,402)&(94,354)&(20,294) \\ \hline
\end{tabular}
\vspace{1mm}

\begin{tabular}{|c|c|c|c|c|c|c|c|c|}\hline $z$&9&10&11&12&13&14&15&16 \\ \hline
&(-396,2198)&(44,438)&(-58,162)&(-4,354)&(76,230)&(6,-158)&(40,590)&(12,22) \\ \hline
\end{tabular}
\vspace{1mm}

\subsection{The Case $D*d$}
\normalsize
This is operator nr. 65 from the list \cite{AESZ}:
\[
{\theta}^{4}-48\,x \left( 6\,\theta+1 \right)  \left( 6\,\theta+5
 \right)  \left( 3\,{\theta}^{2}+3\,\theta+1 \right) +4608\,{x}^{2}
 \left( 6\,\theta+1 \right)  \left( 6\,\theta+5 \right)  \left( 6\,
\theta+7 \right)  \left( 6\,\theta+11 \right)
\]
\scriptsize
$p=3$: \ \ \ \ \ \ \ \ \ \ \ \ \ \ \ \ \ \ \ \ \ \ \ \ \ \ \ \ \ $p=5$:
\vspace{2mm}

$\begin{array}{|c|c|c|}\hline z&1&2 \\ \hline
&(-1,-8)&(-7,16) \\ \hline
\end{array}$\ \ \ \
$\begin{array}{|c|c|c|c|c|}\hline z&1&2&3&4 \\ \hline
&(-26,26)*&(-11,-64)*&(-1,-2)&(14,23) \\ \hline
\end{array}$
\vspace{2mm}

\noindent $p=7$:
\vspace{2mm}

\begin{tabular}{|c|c|c|c|c|c|c|}\hline $z$&1&2&3&4&5&6 \\ \hline
&(-3,-4)&(-12,54)&-&(-9,40)&(-11,66)&(43,-6)* \\ \hline
\end{tabular}
\vspace{2mm}

\noindent $p=11$:
\vspace{2mm}

\begin{tabular}{|c|c|c|c|c|c|c|c|c|}\hline $z$&1&2&3&4&5&6&7&8 \\ \hline
&(67,-538)*&(-19,-58)&(-17,-62)&(-48,106)&(-13,104)'&-&(79,324)&(-62,282) \\ \hline
\end{tabular}
\vspace{1mm}

\begin{tabular}{|c|c|c|}\hline $z$&9&10 \\ \hline
&(-5,-128)&(30,173) \\ \hline
\end{tabular}
\vspace{2mm}

\noindent $p=13$:
\vspace{2mm}

\begin{tabular}{|c|c|c|c|c|c|c|c|c|}\hline $z$&1&2&3&4&5&6&7&8 \\ \hline
&-&(22,178)&-&(24,-70)&(87,428)&(-164,86)*&-&(-35,142) \\ \hline
\end{tabular}
\vspace{1mm}

\begin{tabular}{|c|c|c|c|c|}\hline $z$&9&10&11&12 \\ \hline
&(-58,179)&(47,276)&(33,-86)&(-126,-446)* \\ \hline
\end{tabular}
\vspace{2mm}

\noindent $p=17$:
\vspace{2mm}

\begin{tabular}{|c|c|c|c|c|c|c|c|c|}\hline $z$&1&2&3&4&5&6&7&8 \\ \hline
&(16,30)&(-31,364)&(-23,94)&(40,281)&(22,99)&(59,-24)&(-410,2450)*&(43,592) \\ \hline
\end{tabular}
\vspace{1mm}

\begin{tabular}{|c|c|c|c|c|c|c|c|c|}\hline $z$&9&10&11&12&13&14&15&16 \\ \hline
&(109,472)&(-25,158)&(15,230)&(110,690)&(5,552)&(-198,-1366)*&(-40,342)&(20,-50) \\ \hline
\end{tabular}
\vspace{1mm}

\subsection{The Case $A*f$}
\normalsize
This is operator nr. 133 from the list \cite{AESZ}:
\[
{\theta}^{4}-12\,x \left( 2\,\theta+1 \right) ^{2} \left( 3\,{\theta}^{
2}+3\,\theta+1 \right) +432\,{x}^{2} \left( 2\,\theta+1 \right) ^{2}
 \left( 2\,\theta+3 \right) ^{2}
\]
\scriptsize

$p=3$:\ \ \ \ \ \ \ \ \ \ \ \ \ \ \ \ \ \ \ \ \ \ \ \ \ \ \ \ \ $p=5$:
\vspace{2mm}

$\begin{array}{|c|c|c|}\hline z&1&2 \\ \hline
&(2,10)'&(-1,-2) \\ \hline
\end{array}$\ \ \ \
$\begin{array}{|c|c|c|c|c|}\hline z&1&2&3&4 \\ \hline
&(3,44)&(-6,-6)'&(-3,28)&(-18,42)' \\ \hline
\end{array}$
\vspace{2mm}

\noindent $p=7$:
\vspace{2mm}

\begin{tabular}{|c|c|c|c|c|c|c|}\hline $z$&1&2&3&4&5&6 \\ \hline
&(48,34)*&(9,26)&(1,26)&(17,26)&(64,162)*&(9,26) \\ \hline
\end{tabular}
\vspace{2mm}
\newpage
\noindent $p=11$:
\vspace{2mm}

\begin{tabular}{|c|c|c|c|c|c|c|c|c|}\hline $z$&1&2&3&4&5&6&7&8 \\ \hline
&(-48,210)&-&-&(3,158)&(-36,18)&(-36,82)&(27,70)&(54,266) \\ \hline
\end{tabular}
\vspace{2mm}

\begin{tabular}{|c|c|c|}\hline $z$&9&10 \\ \hline
&(21,-58)&(-54,122) \\ \hline
\end{tabular}
\vspace{1mm}

\noindent $p=13$:
\vspace{2mm}

\begin{tabular}{|c|c|c|c|c|c|c|c|c|}\hline $z$&1&2&3&4&5&6&7&8 \\ \hline
&-&(38,146)&(-47,48)&(-18,-38)&(-192,478)*&(133,660)&(-11,84)&(-34,146) \\ \hline
\end{tabular}
\vspace{1mm}

\begin{tabular}{|c|c|c|c|c|}\hline $z$&9&10&11&12 \\ \hline
&(-18,-166)&(58,242)&(-192,478)*&(50,98) \\ \hline
\end{tabular}
\vspace{2mm}

\noindent $p=17$:
\vspace{2mm}

\begin{tabular}{|c|c|c|c|c|c|c|c|c|}\hline $z$&1&2&3&4&5&6&7&8 \\ \hline
&(48,350)&-&-&(-48,286)&(-9,-260)&-&(72,494)&(-111,524) \\ \hline
\end{tabular}
\vspace{1mm}

\begin{tabular}{|c|c|c|c|c|c|c|c|c|}\hline $z$&9&10&11&12&13&14&15&16 \\ \hline
&(72,622)&(-81,268)&(6,42)&(-48,334)&(42,-54)&(-18,-54)&(-126,570)&- \\ \hline
\end{tabular}
\vspace{1mm}

\subsection{The Case $B*f$}
\normalsize
This is operator nr. 134 from the list \cite{AESZ}:
\[
{\theta}^{4}-9\,x \left( 3\,\theta+1 \right)  \left( 3\,\theta+2
 \right)  \left( 3\,{\theta}^{2}+3\,\theta+1 \right) +243\,{x}^{2}
 \left( 3\,\theta+1 \right)  \left( 3\,\theta+2 \right)  \left( 3\,
\theta+4 \right)  \left( 3\,\theta+5 \right)
\]
\scriptsize
$p=3$:\ \ \ \ \ \ \ \ \ \ \ \ \ \ \ \ \ \ \ \ \ \ \ \ \ \ $p=5$:
\vspace{2mm}

$\begin{array}{|c|c|c|}\hline z&1&2 \\ \hline
&(-4,13)&(5,4) \\ \hline
\end{array}$\ \ \ \
$\begin{array}{|c|c|c|c|c|}\hline z&1&2&3&4 \\ \hline
&(-24,71)&(3,17)&-&(-3,-31) \\ \hline
\end{array}$
\vspace{2mm}

\noindent $p=7$:
\vspace{2mm}

\begin{tabular}{|c|c|c|c|c|c|c|}\hline $z$&1&2&3&4&5&6 \\ \hline
&(11,75)&-&(5,-12)&(-34,-78)*&(-34,-78)*&(5,60)' \\ \hline
\end{tabular}
\vspace{2mm}

\noindent $p=11$:
\vspace{2mm}

\begin{tabular}{|c|c|c|c|c|c|c|c|c|}\hline $z$&1&2&3&4&5&6&7&8 \\ \hline
&(15,218)&(-78,296)&(-12,2)&(-36,194)&(-3,-79)&(69,263)&(-36,113)&- \\ \hline
\end{tabular}
\vspace{1mm}

\begin{tabular}{|c|c|c|}\hline $z$&9&10 \\ \hline
&(-24,107)&(-9,131) \\ \hline
\end{tabular}
\vspace{2mm}

\noindent $p=13$:
\vspace{2mm}

\begin{tabular}{|c|c|c|c|c|c|c|c|c|}\hline $z$&1&2&3&4&5&6&7&8 \\ \hline
&(-1,-171)&(-133,-348)*&(23,114)&(41,159)&(-25,165)&(-109,450)&(-133,-348)*&(32,-48) \\ \hline
\end{tabular}
\vspace{1mm}

\begin{tabular}{|c|c|c|c|c|}\hline $z$&9&10&11&12 \\ \hline
&(98,495)&(-55,99)&(50,33)&(44,306)' \\ \hline
\end{tabular}
\vspace{2mm}

\noindent $p=17$:
\vspace{2mm}

\begin{tabular}{|c|c|c|c|c|c|c|c|c|}\hline $z$&1&2&3&4&5&6&7&8 \\ \hline
&(-12,-322)&(-135,695)&(-105,506)&(-63,227)&(30,434)&(-24,-286)&(45,254)&(-156,857) \\ \hline
\end{tabular}
\vspace{1mm}

\begin{tabular}{|c|c|c|c|c|c|c|c|c|}\hline $z$&9&10&11&12&13&14&15&16 \\ \hline
&(42,92)&(15,-25)&(30,-142)&(12,362)&(-6,236)&(108,641)&(15,461)&(-84,587) \\ \hline
\end{tabular}
\vspace{1mm}

\subsection{The Case $C*f$}
\normalsize
This is operator nr. 135 from the list \cite{AESZ}:
\[
{\theta}^{4}-12\,x \left( 4\,\theta+1 \right)  \left( 4\,\theta+3
 \right)  \left( 3\,{\theta}^{2}+3\,\theta+1 \right) +432\,{x}^{2}
 \left( 4\,\theta+1 \right)  \left( 4\,\theta+3 \right)  \left( 4\,
\theta+5 \right)  \left( 4\,\theta+7 \right)
\]
\scriptsize
$p=3$:\ \ \ \ \ \ \ \ \ \ \ \ \ \ \ \ \ \ \ \ \ \ \ \ \ \ \ \ \ $p=5$:
\vspace{2mm}

$\begin{array}{|c|c|c|}\hline z&1&2 \\ \hline
&(-4,-2)'&(5,10)' \\ \hline
\end{array}$\ \ \ \
$\begin{array}{|c|c|c|c|c|}\hline z&1&2&3&4 \\ \hline
&(-12,22)&(-3,34)&(6,26)&- \\ \hline
\end{array}$
\vspace{2mm}

\noindent $p=7$:
\vspace{2mm}

\begin{tabular}{|c|c|c|c|c|c|c|}\hline $z$&1&2&3&4&5&6 \\ \hline
&(5,-46)&(60,130)*&(52,66)*&(9,58)&(1,62)&- \\ \hline
\end{tabular}
\vspace{2mm}

\noindent $p=11$:
\vspace{2mm}

\begin{tabular}{|c|c|c|c|c|c|c|c|c|}\hline $z$&1&2&3&4&5&6&7&8 \\ \hline
&-&(-24,146)&(-18,38)&(-48,146)&(51,202)&(15,50)&(-24,26)&- \\ \hline
\end{tabular}
\vspace{1mm}

\begin{tabular}{|c|c|c|}\hline $z$&9&10 \\ \hline
&(-78,322)&(-27,-30) \\ \hline
\end{tabular}
\vspace{2mm}

\noindent $p=13$:
\vspace{2mm}

\begin{tabular}{|c|c|c|c|c|c|c|c|c|}\hline $z$&1&2&3&4&5&6&7&8 \\ \hline
&(14,2)&(-54,142)&(44,230)&(11,-124)&(5,-190)&(210,730)*&(-22,-118)&- \\ \hline
\end{tabular}
\vspace{1mm}

\begin{tabular}{|c|c|c|c|c|}\hline $z$&9&10&11&12 \\ \hline
&(44,198)'&(99,436)&(154,-54)&- \\ \hline
\end{tabular}
\vspace{2mm}

\noindent $p=17$:
\vspace{2mm}

\begin{tabular}{|c|c|c|c|c|c|c|c|c|}\hline $z$&1&2&3&4&5&6&7&8 \\ \hline
&-&(-111,688)&(90,494)&(-39,44)&(6,-358)&(42,322)&(-138,810)&- \\ \hline
\end{tabular}
\vspace{1mm}

\begin{tabular}{|c|c|c|c|c|c|c|c|c|}\hline $z$&9&10&11&12&13&14&15&16 \\ \hline
&(-105,412)&(135,698)&-&(6,534)&(-72,622)&(-39,74)&(12,262)&(-36,582)' \\ \hline
\end{tabular}
\vspace{1mm}

\subsection{The Case $D*f$}
\normalsize
This is operator nr. 136 from the list \cite{AESZ}:
\[
{\theta}^{4}-36\,x \left( 6\,\theta+1 \right)  \left( 6\,\theta+5
 \right)  \left( 3\,{\theta}^{2}+3\,\theta+1 \right) +3888\,{x}^{2}
 \left( 6\,\theta+1 \right)  \left( 6\,\theta+5 \right)  \left( 6\,
\theta+7 \right)  \left( 6\,\theta+11 \right)
\]
\scriptsize
$p=3$:\ \ \ \ \ \ \ \ \ \ \ \ \ \ \ \ \ \ \ \ \ \ \ \ \ \ \ $p=5$:
\vspace{2mm}

$\begin{array}{|c|c|c|}\hline z&1&2 \\ \hline
&(-4,13)&(5,4) \\ \hline
\end{array}$ \ \ \ \
$\begin{array}{|c|c|c|c|c|}\hline z&1&2&3&4 \\ \hline
&-&-&(-6,-7)&(-21,67) \\ \hline
\end{array}$
\vspace{2mm}

\noindent $p=7$:
\vspace{2mm}

\begin{tabular}{|c|c|c|c|c|c|c|}\hline $z$&1&2&3&4&5&6 \\ \hline
&(15,-1)&(52,66)*&(5,80)'&-&(9,44)&(60,130)* \\ \hline
\end{tabular}
\vspace{2mm}

\noindent $p=11$:
\vspace{2mm}

\begin{tabular}{|c|c|c|c|c|c|c|c|c|}\hline $z$&1&2&3&4&5&6&7&8 \\ \hline
&(-27,123)&(-9,137)&(-24,62)&(-6,-163)&(-54,246)&(-36,254)&(48,208)&- \\ \hline
\end{tabular}
\vspace{1mm}

\begin{tabular}{|c|c|c|}\hline $z$&9&10 \\ \hline
&(51,126)&(-63,289) \\ \hline
\end{tabular}
\vspace{2mm}

\noindent $p=13$:
\vspace{2mm}

\begin{tabular}{|c|c|c|c|c|c|c|c|c|}\hline $z$&1&2&3&4&5&6&7&8 \\ \hline
&(35,162)&(35,86)&(64,243)&(-60,310)&(-207,688)*&-&-&(20,5) \\ \hline
\end{tabular}
\vspace{1mm}

\begin{tabular}{|c|c|c|c|c|}\hline $z$&9&10&11&12 \\ \hline
&(93,383)&-&(-207,688)*&(5,-259) \\ \hline
\end{tabular}
\vspace{2mm}

\noindent $p=17$:
\vspace{2mm}

\begin{tabular}{|c|c|c|c|c|c|c|c|c|}\hline $z$&1&2&3&4&5&6&7&8 \\ \hline
&(-90,265)&(-99,653)&(6,375)&(-132,580)&(48,230)&(36,394)&(87,335)&(-156,760) \\ \hline
\end{tabular}
\vspace{1mm}

\begin{tabular}{|c|c|c|c|c|c|c|c|c|}\hline $z$&9&10&11&12&13&14&15&16 \\ \hline
&(72,415)&(-45,342)&(-48,74)&(12,-94)&(33,-201)&(15,478)&(9,-225)&(-36,74) \\ \hline
\end{tabular}
\vspace{1mm}

\subsection{The Case $A*g$ }
\normalsize
This is operator nr. 137 from the list \cite{AESZ}:
\[
{\theta}^{4}-4\,x \left( 17\,{\theta}^{2}+17\,\theta+6 \right)
 \left( 2\,\theta+1 \right) ^{2}+1152\,{x}^{2} \left( 2\,\theta+1
 \right) ^{2} \left( 2\,\theta+3 \right) ^{2}
\]
\scriptsize

$p=3$:\ \ \ \ \ \ \ \ \ \ \ \ \ \ \ \ \ \ \ \ $p=5$:
\vspace{2mm}

$\begin{array}{|c|c|c|}\hline z&1&2 \\ \hline
&-&(8,2)* \\ \hline
\end{array}$\ \ \ \
$\begin{array}{|c|c|c|c|c|}\hline z&1&2&3&4 \\ \hline
&-&(-32,62)*&(-6,42)'&(16,-34)* \\ \hline
\end{array}$
\vspace{2mm}

\noindent $p=7$:
\vspace{2mm}

\begin{tabular}{|c|c|c|c|c|c|c|}\hline $z$&1&2&3&4&5&6 \\ \hline
&(6,50)&(80,290)*&(8,2)&(32,-94)*&(16,2)&(6,34) \\ \hline
\end{tabular}
\vspace{2mm}

\noindent $p=11$:
\vspace{2mm}

\begin{tabular}{|c|c|c|c|c|c|c|c|c|}\hline $z$&1&2&3&4&5&6&7&8 \\ \hline
&(-104,-94)*&(-8,98)'&(2,170)&(-64,194)&(-32,2)&(8,2)&-&- \\ \hline
\end{tabular}
\vspace{1mm}

\begin{tabular}{|c|c|c|}\hline $z$&9&10 \\ \hline
&(12,114)'&- \\ \hline
\end{tabular}
\vspace{2mm}

\noindent $p=13$:
\vspace{2mm}

\begin{tabular}{|c|c|c|c|c|c|c|c|c|}\hline $z$&1&2&3&4&5&6&7&8 \\ \hline
&(-108,-698)*&(14,146)&-&(-56,174)'&-&(-160,30)*&(36,278)&(36,118) \\ \hline
\end{tabular}
\vspace{1mm}

\begin{tabular}{|c|c|c|c|c|}\hline $z$&9&10&11&12 \\ \hline
&(66,322)&(-36,22)&(16,-114)&(24,206) \\ \hline
\end{tabular}
\vspace{2mm}

\noindent $p=17$:
\vspace{2mm}

\begin{tabular}{|c|c|c|c|c|c|c|c|c|}\hline $z$&1&2&3&4&5&6&7&8 \\ \hline
&(-88,494)'&(-356,1478)*&(-40,14)&(92,326)'&(4,-154)&(88,350)&(10,-430)'&(6,-174) \\ \hline
\end{tabular}
\vspace{1mm}

\begin{tabular}{|c|c|c|c|c|c|c|c|c|}\hline $z$&9&10&11&12&13&14&15&16 \\ \hline
&(6,210)&(-148,854)'&-&(56,206)&(-92,566)&(-182,1010)&(224,-898)*&(64,62) \\ \hline
\end{tabular}
\vspace{1mm}

\subsection{The Case $B*g$ }
\normalsize
This is operator nr. 138 from the list \cite{AESZ}:
\[
{\theta}^{4}-3\,x \left( 3\,\theta+1 \right)  \left( 3\,\theta+2
 \right)  \left( 17\,{\theta}^{2}+17\,\theta+6 \right) +648\,{x}^{2}
 \left( 3\,\theta+1 \right)  \left( 3\,\theta+2 \right)  \left( 3\,
\theta+4 \right)  \left( 3\,\theta+5 \right)
\]
\scriptsize

$p=3$:\ \ \ \ \ \ \ \ \ \ \ \ \ \ \ \ \ \ \ \ \ \ \ \ \ \ \ \ \ \ $p=5$:
\vspace{2mm}

$\begin{array}{|c|c|c|}\hline z&1&2 \\ \hline
&(-4,-2)'&(5,10)' \\ \hline
\end{array}$\ \ \ \
$\begin{array}{|c|c|c|c|c|}\hline z&1&2&3&4 \\ \hline
&(18,-22)*&(-33,68)*&(-9,14)&- \\ \hline
\end{array}$
\vspace{2mm}

\noindent $p=7$:
\vspace{2mm}

\begin{tabular}{|c|c|c|c|c|c|c|}\hline $z$&1&2&3&4&5&6 \\ \hline
&(5,-66)&(32,96)&(-46,18)*&(23,96)&(-13,-12)&- \\ \hline
\end{tabular}
\vspace{2mm}

\noindent $p=11$:
\vspace{2mm}

\begin{tabular}{|c|c|c|c|c|c|c|c|c|}\hline $z$&1&2&3&4&5&6&7&8 \\ \hline
&(-120,98)*&(6,-37)&(-24,89)&(60,206)'&-&-&-&(72,-478)* \\ \hline
\end{tabular}
\vspace{1mm}

\begin{tabular}{|c|c|c|}\hline $z$&9&10 \\ \hline
&(-9,-10)&(-39,134) \\ \hline
\end{tabular}
\vspace{2mm}
\newpage
\noindent $p=13$:
\vspace{2mm}

\begin{tabular}{|c|c|c|c|c|c|c|c|c|}\hline $z$&1&2&3&4&5&6&7&8 \\ \hline
&(-31,6)'&-&(-190,450)*&(86,321)&(-103,-768)*&(-4,222)&(-1,288)'&(-16,66)' \\ \hline
\end{tabular}
\vspace{1mm}

\begin{tabular}{|c|c|c|c|c|}\hline $z$&9&10&11&12 \\ \hline
&(-16,210)&-&(14,33)&(41,294)' \\ \hline
\end{tabular}
\vspace{2mm}

\noindent $p=17$:
\vspace{2mm}

\begin{tabular}{|c|c|c|c|c|c|c|c|c|}\hline $z$&1&2&3&4&5&6&7&8 \\ \hline
&(-18,506)&(45,344)&(63,146)&(36,83)&(-171,902)&-&(-432,2846)*&(-150,812) \\ \hline
\end{tabular}
\vspace{1mm}

\begin{tabular}{|c|c|c|c|c|c|c|c|c|}\hline $z$&9&10&11&12&13&14&15&16 \\ \hline
&(-66,164)&(414,2522)*&(-15,-160)&(-57,146)&(-36,-241)&(3,524)&-&(-6,182) \\ \hline
\end{tabular}
\vspace{1mm}

\subsection{The Case $C*g$}
\normalsize
This is operator nr. 139 from the list \cite{AESZ}:
\[
{\theta}^{4}-4\,x \left( 4\,\theta+1 \right)  \left( 4\,\theta+3
 \right)  \left( 17\,{\theta}^{2}+17\,\theta+6 \right) +1152\,{x}^{2}
 \left( 4\,\theta+1 \right)  \left( 4\,\theta+3 \right)  \left( 4\,
\theta+5 \right)  \left( 4\,\theta+7 \right)
\]
\scriptsize
$p=3$:\ \ \ \ \ \ \ \ \ \ \ \ \ \ \ \ \ \ \ \ \ \ \ \ \ \ \ \ \ \ \ \ \ $p=5$:
\vspace{2mm}

$\begin{array}{|c|c|c|}\hline z&1&2 \\ \hline
&(-4,10)&(-4,-14)*\\ \hline
\end{array}$\ \ \ \
$\begin{array}{|c|c|c|c|c|}\hline z&1&2&3&4 \\ \hline
&(-28,38)*&-&(18,-22)*&(-4,14) \\ \hline
\end{array}$
\vspace{2mm}

\noindent $p=7$:
\vspace{2mm}

\begin{tabular}{|c|c|c|c|c|c|c|}\hline $z$&1&2&3&4&5&6 \\ \hline
&(88,354)*&(2,-46)&-&(68,194)*&(-18,22)&(-6,58)' \\ \hline
\end{tabular}
\vspace{2mm}

\noindent $p=11$:
\vspace{2mm}

\begin{tabular}{|c|c|c|c|c|c|c|c|c|}\hline $z$&1&2&3&4&5&6&7&8 \\ \hline
&(-14,194)&(-140,338)*&(72,-478)*&(50,290)&(-8,130)&(-58,202)&(10,-70)&(-24,-6) \\ \hline
\end{tabular}
\vspace{1mm}

\begin{tabular}{|c|c|c|}\hline $z$&9&10 \\ \hline
&(16,106)&(-24,2) \\ \hline
\end{tabular}
\vspace{2mm}

\noindent $p=13$:
\vspace{2mm}

\begin{tabular}{|c|c|c|c|c|c|c|c|c|}\hline $z$&1&2&3&4&5&6&7&8 \\ \hline
&(-20,294)&(2,-126)&-&-&-&(32,134)&(38,174)&(202,618)* \\ \hline
\end{tabular}
\vspace{1mm}

\begin{tabular}{|c|c|c|c|c|}\hline $z$&9&10&11&12 \\ \hline
&(30,-62)&(224,926)*&(-22,38)&(20,270) \\ \hline
\end{tabular}
\vspace{2mm}

\noindent $p=17$:
\vspace{2mm}

\begin{tabular}{|c|c|c|c|c|c|c|c|c|}\hline $z$&1&2&3&4&5&6&7&8 \\ \hline
&-&(-22,338)&(-128,746)&(-44,86)&(-50,74)&(-44,74)&(-52,14)&(316,758)* \\ \hline
\end{tabular}
\vspace{1mm}

\begin{tabular}{|c|c|c|c|c|c|c|c|c|}\hline $z$&9&10&11&12&13&14&15&16 \\ \hline
&(-208,1186)*&(40,110)&(-22,-94)&(164,818)&(8,-370)&(-52,218)&(-182,1010)&(-64,302) \\ \hline
\end{tabular}
\vspace{1mm}

\subsection{The Case $D*g$}
\normalsize
This is operator nr. 140 from the list \cite{AESZ}:
\[
{\theta}^{4}-12\,x \left( 6\,\theta+1 \right)  \left( 6\,\theta+5
 \right)  \left( 17\,{\theta}^{2}+17\,\theta+6 \right) +10368\,{x}^{2}
 \left( 6\,\theta+1 \right)  \left( 6\,\theta+5 \right)  \left( 6\,
\theta+7 \right)  \left( 6\,\theta+11 \right)
\]
\scriptsize
$p=3$:\ \ \ \ \ \ \ \ \ \ \ \ \ \ \ \ \ \ \ \ \ \ \ \ \ \ \ \ \ \ $p=5$:
\vspace{2mm}

$\begin{array}{|c|c|c|}\hline z&1&2 \\ \hline
&(-4,-2)'&(5,10)' \\ \hline
\end{array}$\ \ \ \
$\begin{array}{|c|c|c|c|c|}\hline z&1&2&3&4 \\ \hline
&(-26,26)*&(19,-16)*&-&- \\ \hline
\end{array}$
\vspace{2mm}

\newpage
\noindent $p=7$:
\vspace{2mm}

\begin{tabular}{|c|c|c|c|c|c|c|}\hline $z$&1&2&3&4&5&6 \\ \hline
&(24,64)&-&(53,74)*&(29,86)&(26,-142)*&(9,-4) \\ \hline
\end{tabular}
\vspace{2mm}

\noindent $p=11$:
\vspace{5mm}

\begin{tabular}{|c|c|c|c|c|c|c|c|c|}\hline $z$&1&2&3&4&5&6&7&8 \\ \hline
&(13,50)&(-17,2)&(8,-58)&-&(14,83)&(160,578)*&(-68,257)&(-32,172) \\ \hline
\end{tabular}
\vspace{1mm}

\begin{tabular}{|c|c|c|}\hline $z$&9&10 \\ \hline
&(-128,194)*&(-37,38) \\ \hline
\end{tabular}
\vspace{2mm}

\noindent $p=13$:
\vspace{2mm}

\begin{tabular}{|c|c|c|c|c|c|c|c|c|}\hline $z$&1&2&3&4&5&6&7&8 \\ \hline
&(-198,562)*&(-4,-202)&(-4,218)&(49,294)&-&(-193,492)*&(84,386)&(24,163) \\ \hline
\end{tabular}
\vspace{1mm}

\begin{tabular}{|c|c|c|c|c|}\hline $z$&9&10&11&12 \\ \hline
&(-9,214)&(24,211)&(19,-36)&(54,170) \\ \hline
\end{tabular}
\vspace{2mm}

\noindent $p=17$:
\vspace{2mm}

\begin{tabular}{|c|c|c|c|c|c|c|c|c|}\hline $z$&1&2&3&4&5&6&7&8 \\ \hline
&(-10,-394)&(71,368)'&(157,908)'&(112,431)&(-47,122)&(3,-144)&(-350,1370)*&(-154,716) \\ \hline
\end{tabular}
\vspace{1mm}

\begin{tabular}{|c|c|c|c|c|c|c|c|c|}\hline $z$&9&10&11&12&13&14&15&16 \\ \hline
&(-38,-76)&(236,-682)*&(-63,-24)&-&-&(5,-394)&(49,320)&(-38,338) \\ \hline
\end{tabular}

KIRA SAMOL: Fachbereich Mathematik 17, AG Algebraische Geometrie, Johannes Gutenberg-Universit\"at, D-55099 Mainz, GERMANY;\\
{\tt samol@uni-mainz.de}\\

DUCO VAN STRATEN:  Fachbereich Mathematik 17, AG Algebraische Geometrie, Johannes Gutenberg-Universit\"at, D-55099 Mainz, GERMANY;\\
{\tt straten@mathematik.uni-mainz.de}

\end{document}